\documentclass[12pt]{article}
\usepackage[dvips]{graphics}
\usepackage[pdftex]{graphicx}
\usepackage{graphicx}
\usepackage{color}
\usepackage{epsfig}
\usepackage{amssymb}
\usepackage{amsmath}
\usepackage{color}

\DeclareMathSizes{12}{12}{10}{10}

\newtheorem{theorem}{Theorem}

\newtheorem{lemma}[theorem]{Lemma}

\newtheorem{remark}[theorem]{Remark}

\newcommand{\kom}[1]{}
\renewcommand{\kom}[1]{{\bf [#1]}}
\definecolor{blau}{rgb}{0.1,0.0,0.9}

\newcounter{komcounter}
\numberwithin{komcounter}{section}

\begin{document}

\title{{ Estimates of size of cycle in a predator-prey system II}}

\author{Niklas L.P. Lundstr\"{o}m$^{\tiny 1}$ and
Gunnar J. S\"oderbacka$^{\tiny 2}$
 \linebreak \\\\
$^{\tiny 1}$\it \small Department of Mathematics and Mathematical Statistics, Ume{\aa} University\\
\it \small SE-90187 Ume{\aa}, Sweden\/{\rm ;}
\it \small niklas.lundstrom@umu.se  \linebreak \\\\
$^{\tiny 2}$\it \small {\AA}bo Akademi, 20500 {\AA}bo, Finland\/{\rm ;}
\it \small gsoderba@abo.fi \linebreak \\\\}


\date{}


\maketitle

\begin{abstract}
We prove estimates for the maximal and minimal predator and prey populations on the unique limit cycle in a standard predator-prey system.
Our estimates are valid when the cycle exhibits small predator and prey abundances and large amplitudes.
The proofs consist of constructions of several Lyapunov-type functions and derivation of a large number of non-trivial estimates,
and should be of independent interest.
This study generalizes results proved by the authors in \cite{ns} to a wider class of systems and, in addition, it gives simpler proofs of some already known estimates.
\\


\end{abstract}

\section{Introduction}

%
%


We consider the predator-prey system
\begin{align}\label{eqmain}
&\frac{ds}{d\tau} = (h(s)-x)s, \notag\\
&\frac{dx}{d\tau} = m (s-\lambda) x, \quad \text{where} \quad h(s)=(1-s)(s+a), \tag{$\star$} \\
&\text{assuming} \quad
a,\lambda,m\in \mathbf{R}_+\quad \text{and} \quad 2 \lambda + a < 1,\notag
\end{align}
%


\noindent
in which $s = s(\tau)$ and $x = x(\tau)$ denote the prey and predator, respectively.
Systems of this type were first introduced in \cite{RM63} as a more realistic modification of the original Lotka-Volterra system.
Indeed, \eqref{eqmain} is equivalent with the following
Rosenzweig-MacArthur system:
\begin{align}\label{Rosenzweig-MacArthur}
 \frac{dS}{dt} &\,=\, r S \left(1-\frac{S}{K}\right) -  \frac{q X S}{H + S}, \notag\\
 \frac{dX}{dt} &\,=\,  \frac{p X S}{H + S} - d X.
\end{align}
Here, $S \,=\, S\left(t\right)$ and $X \,=\, X\left(t\right)$ denote the population densities of prey and predator,
and $r, K, q, H, p$ and $d$ are positive parameters.
The biological meanings of the parameters are the following:
$r$ is the intrinsic growth rate of the prey;
$K$ is the prey carrying capacity;
$q$ is the maximal consumption rate of predators;
$H$ is the amount of prey needed to achieve one-half of $q$;
$p$ is the efficiency with which predators convert consumed prey into new predators;
and $d$ is the per capita death rate of predators.

System \eqref{Rosenzweig-MacArthur} transforms into system \eqref{eqmain}
when introducing the non-dimensional quantities
\begin{align}\label{def-parameters}
\tau \,=\, \int \frac{r K}{H + S\left(t\right)} dt, \quad s &\,=\, \frac{S}{K}, \quad x \,=\, \frac{q X}{r K}, \quad a \,=\, \frac{H}{K},  \notag\\
m \,=\, \frac{p - d}{r} \quad &\text{and} \quad \lambda \,=\, \frac{d H}{(p-d) K}.
\end{align}
Rosenzweig-MacArthur-type systems have been widely used in ecological applications and have been frequently studied by mathematicians,
see e.g. \cite{C81, HS09, ns} and the references therein.
%
%
System \eqref{eqmain} always has a unique positive equilibrium at
$\left(x,s\right) = \left(\left(1-\lambda\right)\left(\lambda + a\right), \lambda\right)$
which attracts the whole positive space when $2\lambda + a \,>\, 1$.
At $2\lambda + a \,=\, 1$ there is a Hopf bifurcation in which the equilibrium loses stability and a stable limit cycle, surrounding the equilibrium, is created.
In particular, for $2\lambda + a \,<\, 1$ the equilibrium is a source and the cycle attracts the whole positive space (except the source) \cite{C81}.
For results of uniqueness of limit cycles in more general system we refer to \cite{H00, H88, HM89, KF88}.

\begin{figure}[ht!]
\begin{center}
\includegraphics[width=8.0cm,height=7.0cm]{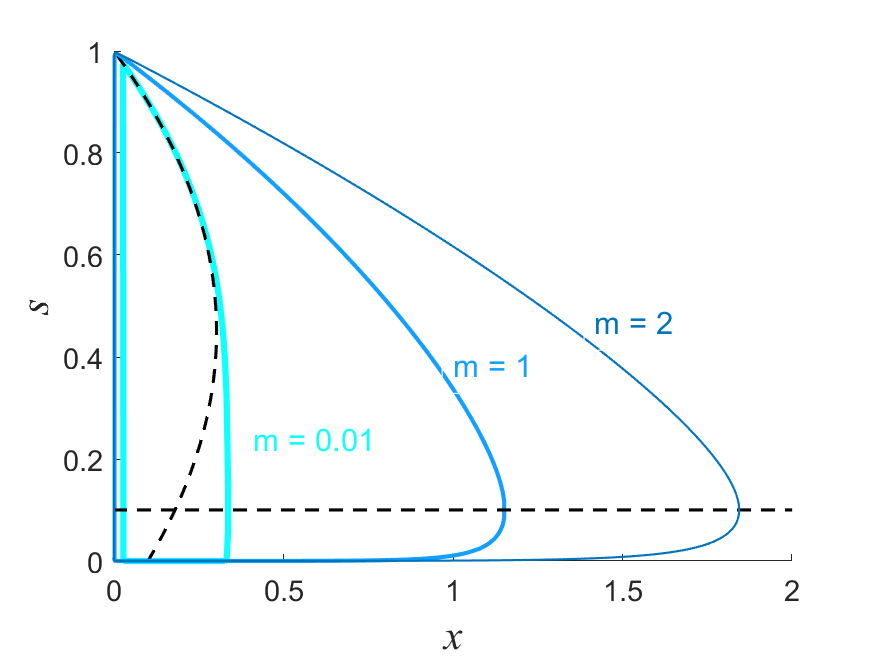}
\includegraphics[width=7.0cm,height=6.5cm]{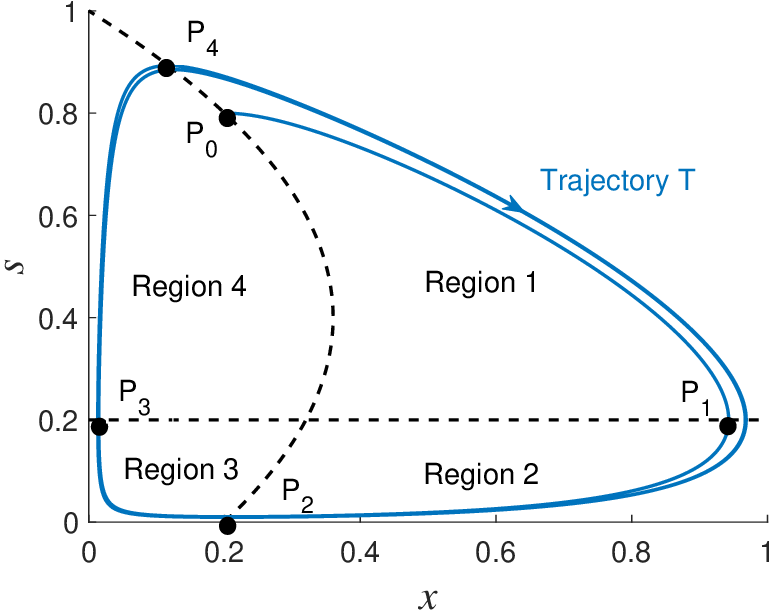}
\end{center}
\caption{Left panel: The unique limit cycle of system \eqref{eqmain} for $a = \lambda = 1/10$ and $m = 0.01$ (thick, light), $m = 1$ and $m = 2$ (thin, dark).
Right panel: Trajectory $T$, the definitions of Regions 1 - 4 and points $P_1$ - $P_4$.
Isoclines $x = h(s)$ and $s = \lambda$ are marked with black dashed curves.}
\label{fig:traj}
\end{figure}

In this paper we prove estimates for this unique stable limit cycle when $a$ and $\lambda$ are small,
and our main contribution is to carry forward the technics from \cite{ns}, valid only for $m =1$, to hold for any $m \in \mathbf{R}_+$.
Figure \ref{fig:traj} shows the unique limit cycle of system \eqref{eqmain} for different values of $m$ when $a = \lambda = 1/10$.
To complete our estimates for the whole cycle we will assume that either
\begin{align}\label{eq:assumption_parameters}
a \leq \frac{1}{20}, \quad \lambda \leq \frac{1}{20} \quad \text{or that} \quad a \leq \frac{1}{10}, \quad \lambda \leq \frac{1}{100}.  \tag{$\star\star$}
\end{align}
We remark that some results for greater $a$ and $\lambda$ can be derived from general results on estimates of limit cycles in \cite{Angelis75}.
To state our main theorem, define
\begin{align*}
\check x &=  \max_{z \in [\frac{1-a}{2},0.8]} h(z) + m\, (z -\lambda\, (1-\ln \lambda +\ln z )),\\
\hat{x} &= \frac{(1+m+a-m\,\lambda)(1+2m+a-2m\, \lambda)}{2(m+1)+a-m\, \lambda},
\end{align*}
and $\check z_1 = z_1(\hat x/a)$, $\hat z_2 = z_2(\check x/h(\lambda))$ in which, for $i = 1,2,$
\begin{align*}
z_i(y) &=\frac{1-d_i\, Y-\sqrt{(1-d_i\, Y)^2 - 4 c_i\, Y}}{2c_i\, Y},\quad \text{where}\quad Y = y\, e^{-y},\\
d_i&=e-1-c_i\, e, \quad  c_1=\frac{e-2}{e-1} \quad \text{and} \quad  c_2=\frac{1}{e}.
\end{align*}
We remark that $\check z_1, \hat z_2 \in(1, e)$ and that $z_i(y)$, $i = 1,2$, decreases towards 1 with increasing $y > 1$.
Our main result is the following theorem: \\

\noindent
{\bf Theorem A}
\emph{ Let $x_{max}$ and $s_{max}$ be the maximal $x$- and $s$-values and let $x_{min}$ and $s_{min}$ be the
minimal $x$- and $s$-values on the unique limit cycle of system \eqref{eqmain} under assumption \eqref{eq:assumption_parameters}.
Then the predator satisfies
\begin{align*}
 \check{x} \, < \, x_{\text{max}} \, < \, \hat{x} \quad\text{and}\quad
  \check z_1 \, \hat x \, e^{ -\frac{\hat x}{a}} \,<\, x_\text{min} \,<\, \hat z_2 \, \check x \, e^{-\frac{\check x}{h(\lambda)}},
\end{align*}
and the prey satisfies
\begin{align*}
0.8 \, < \, s_{\text{max}} \, < \, 1 \quad\text{and}\quad \lambda \,e^{- \frac{1}{m \lambda}\left(\hat x - a - a\ln\frac{\hat x}{a}\right) - 1 } \,< \, s_\text{min} \, < \, \lambda \, e^{- \frac{1}{m \lambda}\left(\check x - a - h(\lambda)\ln\frac{\check x}{a}\right) }.
\end{align*}
%
}

The proof of Theorem A is given in Sections \ref{sec:reg1}, \ref{sec:lower-estimates} and \ref{sec:reg4}.
In different parts of the state-space, we state and prove several results which are valid for wider ranges of parameter values than those needed for Theorem A.
In the case $m=1$ sharper estimates are given in
\cite{ns}, and related estimates are given in \cite{HS09}.
We remark that in Region 2 and Region 3 (see definitions below and Fig. \ref{fig:traj}) the present paper uses a much simpler technique than the corresponding one used by the authors in \cite{ns}. We here adopt some rather simple analytical approximations of Lotka-Volterra integrals,
see \cite{LS23} and Lemma \ref{le:julelemma}.

Biologically, our main results imply that the modeled population exhibits very small predator and prey abundances
during a long portion in time of the cycle,
indicating that the population suffers a relatively high risk of going extinct because of random perturbations.
That the modelled predator and prey populations survives even though reaching extremely low densities may be unrealistic,
but it should be noted that in reality, predators may survive even though prey biomass is very small by feeding on alternative resources.
From \eqref{def-parameters} we see that an increase in the carrying capacity $K$
(when other parameters are fixed) implies a decrease in both parameters $a$ and $\lambda$.
Thus, if $K$ is large enough then assumption \eqref{eq:assumption_parameters} will be satisfied.
Our results therefore show that the population becomes small,
and hence vulnerable,
when the carrying capacity $K$ becomes large.
This result is in line with the paradox of enrichment \cite{R71}.
Indeed, Rosenzweig argues that enrichment of the environment (larger carrying capacity $K$)
leads to destabilization.
%

The Rosenzweig-MacArthur system is a very simplified model of reality,
but the results of the present study may be useful also when investigating dynamics in more complex and less simplified models,
such as systems modeling seasonally dependence, see e.g. \cite{BHW83,EOS96,RMK93}.
Theorem A may also useful when studying systems of type
many predators - one prey like in \cite{LS23,osipovAlushta,osipovIJBC,osipoveuler,S23}.
Often such systems have limit cycles and the stability or instability of these cycles are important for understanding
the extinction or coexistence possibilities for the predators.
It occurs that the behaviour for small prey biomass on the cycles can play an important role in determining the stability,
see \cite{LS23} for detailed information and \cite{S23} where main ideas are presented.
Theorem A gives such estimates in case only one predator survives.


To prove Theorem A we have constructed several Lyapunov-type functions and derived a large number of non-trivial estimates.
We believe that these methods and constructions have values also beyond this paper
as they present methods and ideas that, potentially, can be useful for proving analogous results
for dynamics in similar systems as well as in more complex systems.

\subsection*{The case of $m$ small -- A canard-type solution}




The case of $m$ small is a suitable introductory case before going into the proof of Theorem A.
Therefore, to get the reader familiar with the ideas for the kinds of estimates which we construct in detail later on we derive an approximation of the
limit cycle as $m \to 0$, which turns the dynamic into a canard-type limit cycle.
With $m$ small, the motion is always slow in the $x$-direction, forcing the cycle to either follow closely the isocline $x = h(s)$ with slow motion,
or jump nearly vertically back and forth towards $s = 0$.

Consider a trajectory $M$ starting at the point $\left(\frac{(1+a)^2}{4}, \frac{1-a}{2}\right)$ which is the top of the parabola $x = h(s)$,
forming the isocline $s' = 0$.
As $m$ approaches zero, the corresponding trajectory $M$ will approach a vertical dip from this point towards the $x$-axis at $\left(x_\epsilon, \epsilon\right)$, with $\epsilon \to 0$, $x_\epsilon \to \frac{(1+a)^2}{4} \to x_{max}$ as $m \to 0$, see the trajectory for $m = 0.01$ in Figure \ref{fig:traj}.

To analyze the trajectory in the region $s < \epsilon$ we take off from the phase-plane equation
\begin{align}\label{eq:phase}
\frac{ds}{dx} = \frac{(h(s) - x)s}{m x (s - \lambda)},
\end{align}
which, through the approximation $h(s) \approx a$ valid for $s < \epsilon$, simplifies to a separable ODE.
%
%
Indeed, let
\begin{align*}
U(x,s) = \frac{1}{m} (x - a \ln x) + s - \lambda \ln s
\end{align*}
and observe that trajectory $T$ approximately satisfies
\begin{align*}
U(x, s(x)) \approx U\left(x_{max},\epsilon\right)
\end{align*}
in the region $s < \epsilon$.
Therefore, when $T$ intersects its second time $s = \epsilon$ at approximately its minimal $x$-value at $(x_{min}, \epsilon)$ we have
$
U(x_{min}, \epsilon) \approx U\left(x_{max},\epsilon\right),
$
i.e.,
\begin{align*}
x_{min} - a \ln x_{min}   \approx x_{max} - a \ln x_{max}.
\end{align*}
Neglecting the first term on the left hand side yields
\begin{align}\label{eq:msmall1}
x_{min} \approx x_{max} e^{ -\frac{x_{max}}{a}} \quad\text{in which}\quad x_{max} \approx \frac{(1+a)^2}{4}.
\end{align}
Finally, from $(x_{min}, \lambda)$ trajectory $M$ will, in the limit $m \to 0$, jump up towards the isocline $x = h(s)$ near the point
$(x_{min}, s_{max})$ where
\begin{align}\label{eq:msmall2}
s_{max} \approx \frac{1-a}{2} + \sqrt{\frac{(1-a)^2}{4} - a + x_{min}}.
\end{align}
From this point we move close to the isocline $x = h(s)$ until the top of $h(s)$ to reach the initial condition of $M$.

We may also derive an estimate for $s_{min}$ at the intersection of $M$ with the isocline $x = h(s)$ near $(a,0)$ by solving
\begin{align*}
U(a, s_{min}) = U\left(x_{max},\epsilon\right),
\end{align*}
which gives,
\begin{align}\label{eq:msmall3}
s_{min} \approx 
e^{\frac{v - a(\ln a - 1)}{m\lambda}} \quad \text{in which}\quad v = a\ln x_{max} - x_{max}.
\end{align}

We remark that estimates \eqref{eq:msmall1}--\eqref{eq:msmall3} can easily be generalized to hold for trajectories of any systems on the form \eqref{eqmain}
 but when allowing for $h(s)$ to satisfy only the condition $$h(0) > 0 = h(1), \quad h''(s) < 0,$$
introduced by Osipov in \cite{osipovAlushta} and \cite{osipovIJBC}.
We also observe that here, and for small lambda in general,
we have slow-fast motions like in \cite{HW20} with the difference that we often also have a slow motion near the equilibrium $(x,s) = (0,1)$.

We now turn to the proof of Theorem A, giving upper and lower estimates for $x_{max}$, $s_{min}$, $x_{min}$ and $s_{max}$ for all $m \in \mathbf{R}_+$.
%
%
%

\subsection*{The proof of Theorem A}

To outline the proof of Theorem A we first observe that the coordinate axes are invariant,
and hence the region $x,s > 0$ is also invariant. Therefore, we consider solutions only for
positive $x$ and $s$. Moreover, system \eqref{eqmain} has isoclines at $x = h (s)$ and $s = \lambda$, which lead
us to split the proof by introducing the following four regions:
\begin{itemize}
\item [] Region 1, where $x > h (s)$, $s > \lambda$ and $x$ is growing and $s$ decreasing.
\item [] Region 2, where $x > h (s)$, $s < \lambda$ and both $x$ and $s$ decrease.
\item [] Region 3, where $x < h (s)$, $s < \lambda$ and $x$ decreases and $s$ grows.
\item [] Region 4, where $x < h (s)$, $s > \lambda$ and both $x$ and $s$ increase.
\end{itemize}
Any trajectory starting in Region 1 will enter Region 2 from where it will enter Region 3
and then Region 4 and finally Region 1 again, and the behaviour repeats infinitely.
Figure \ref{fig:traj} illustrates the four regions together with isoclines and points which will be used in the proof
of Theorem A. Behaviour and estimates for trajectories in different regions are examined in
different subsections.

Theorem A follows from four statements, Statement 1--4, which we prove in Sections \ref{sec:reg1}--\ref{sec:reg4} below.
Each statement gives estimates for a trajectory of system \eqref{eqmain} assuming \eqref{eq:assumption_parameters} which we denote by $T$.
Trajectory $T$ stars at $P_0 = \left(h(s_0),s_0\right)$ on the isocline $x = h(s)$ between Region 1 and Region 4.
It will intersect isocline $s = \lambda$ at $P_1 = (x_1, \lambda)$,
isocline $x = h(s)$ at $P_2 = (h(s_2), s_2)$,
and then again isocline $s = \lambda$ at $P_3 = (x_3, \lambda)$.
Finally, trajectory $T$ returns to Region 1 from Region 4 at another intersection with isocline $x = h(s)$ at $P_4 = (h(s_4),s_4)$, see Figure \ref{fig:traj}.\\

\noindent
\emph{Proof of Theorem A.}
Assumption \eqref{eq:assumption_parameters} ensures the validity of Statements 1--4.
Statement 4 implies that $s_4 > 0.8$ when $s_0 = 0.8$ and thus it must follow that the unique limit cycle of system \eqref{eqmain} can
be estimated by Statements 1--3 with $s_0 = 0.8$.
Letting $\check x$ denote $\check x_1$ when $s_0 = 0.8$ and $\hat x = \hat x_1$
(where $\check x_1$ and $\hat x_1$ are defined in Statement 1)
finishes the proof.
$\hfill\Box$\\



\section{Estimates in Region 1}

\label{sec:reg1}

\setcounter{theorem}{0}
\setcounter{equation}{0}


\begin{lemma}\label{th:uppis}
At the intersection of trajectory $T$ with the isocline $s = \lambda$ at $P_1 = (x_1,\lambda)$ it holds that
\begin{equation*}
x_1 < \frac{(1+m+a-m\,\lambda)(1+2m+a-2m\, \lambda)}
{2(m+1)+a-m\, \lambda}.
\end{equation*}
%
%
\end{lemma}


\noindent
{\it Proof}. Let $v=1-s$ and consider the function
\begin{equation*}
V = x - \frac{A\, v}{1+B\, v} \quad \text{with} \quad A=1+m+a-m\,\lambda \quad \text{and}\quad B=\frac{1+m\,\lambda}{1+a+2m\, (1-\lambda)}.
\end{equation*}
The choice of $A$ is natural as it forces $x = A v$ in the direction of the escaping eigenvector at the saddle (0,1).
The derivative $V'$ with respect to time evaluated at $V=0$ yields
\begin{equation*}
V'= \frac{A\, v^3\, (C_1\, v+C_0)}{b^2 \, (1+Bv)^3} 
\end{equation*}
where $b=1+a+2m\, (1-\lambda)$ (the denominator of $B$) and
\begin{equation*}
C_1 + C_0 = -m\, \lambda\, (a+2+2m-m\, \lambda)^2
\end{equation*}
and
\begin{equation*}
C_0 = -m^3\, \lambda (\lambda^2 -5\lambda + 4) +
m^2\, \lambda ((2a+6)\lambda -8-4a) +c
\end{equation*}
where
\begin{equation*}
c= - m\, \lambda \,( 3 + 5a + a^2) -1-m-a .
\end{equation*}
We observe that $C_0<0$ for $0\leq \lambda <1$ and that
$C_0$ is the value of $C_1\, v+C_0$ for $v=0$ and $C_1 + C_0$
is the value of $C_1\, v+C_0$ for $v=1$.
We conclude $V'<0$
from which we get Lemma \ref{th:uppis} by substituting $v=1$ into
$\frac{A\, v}{1+B\, v}$. $\hfill\Box$\\

\begin{remark}\label{re:better}
Substituting $v=1-\lambda$ into $\frac{A\, v}{1+B\, v}$
we get the somewhat better estimate
\begin{equation*}
x_1 < \frac{(1+a+2m\, L)(1+a+ m \, L)\, L}
{1+a+(1+2m+m\lambda)\,L}
\end{equation*}
where $L=1-\lambda$.
Further if $a,\lambda \to 0$ we get the simple estimate $x_1 < m+0.5$.
\end{remark}


\begin{lemma}\label{le:lemma2}
At the intersection of trajectory $T$ with the isocline $s = \lambda$ at $P_1 = (x_1,\lambda)$ it holds that
\begin{equation}\label{xmax2}
x_1 > x_0 +m\, (s_0 -\lambda\, (1-\ln \lambda +\ln s_0 )).
\end{equation}
\end{lemma}

\noindent
{\it Proof}. Consider
$$
V=m\, (s-\lambda \ln s) + x.
$$
Differentiating with respect to time gives $V'=m\, (s-\lambda ) h(s) >0$.
Thus $V > V_0$ on the trajectory, if
$V_0$ is the value of $V$ at $(x_0,s_0)$.
From this we immediately get (\ref{xmax2}). $\hfill \Box$

\begin{remark}
For $s_0=1,\, x_0=0$ we get
$x_1 > m\, (1-\lambda (1-\ln\lambda ))$, which is an
estimate for the maximal $x$-value of the unstable manifold
of the equilibrium $(0,1)$.
\end{remark}

%

We remark that our lower estimate in \eqref{xmax2} is good for large $m$ but it may drop below the trivial lower estimate $\frac{(1+a)^2}{4}$ when $m$ is small.
However, as long as $s_0 > \frac{1-a}{2}$ a possibly better lower barrier for $T$ is produced by replacing $s_0$ with a smaller value in \eqref{xmax2},
i.e. starting the barrier closer to the tip of the parabola $x = h(s)$.
Indeed, for $s_0 > \frac{1-a}{2}$ a better lower bound $\check x_1$ for $x_1$ is given by
$$
\check x_1 = \max_{z \in [\frac{1-a}{2},s_0]} h(z) + m\, (z -\lambda\, (1-\ln \lambda +\ln z )).
$$

We next summarize the above results implication for trajectory $T$ in a statement.
Indeed, the below estimate follows immediately from Lemma \ref{th:uppis} and Lemma \ref{le:lemma2}.\\

\noindent
{\bf Statement 1}
\emph{
At the intersection of trajectory $T$ with the isocline $s = \lambda$ at $P_1 = (x_1,\lambda)$ it holds that}
\begin{equation*}
\check x_1 < x_1 < \hat x_1
\end{equation*}
\emph{where} 
$$
\check x_1 =  \max_{z \in [\frac{1-a}{2},s_0]} h(z) + m\, (z -\lambda\, (1-\ln \lambda +\ln z ))
$$
\emph{and}
$$
\hat x_1 = \frac{(1+m+a-m\,\lambda)(1+2m+a-2m\, \lambda)}
{2(m+1)+a-m\, \lambda}.
$$

\begin{figure}[h]
\begin{center}
\includegraphics[scale = 0.5]{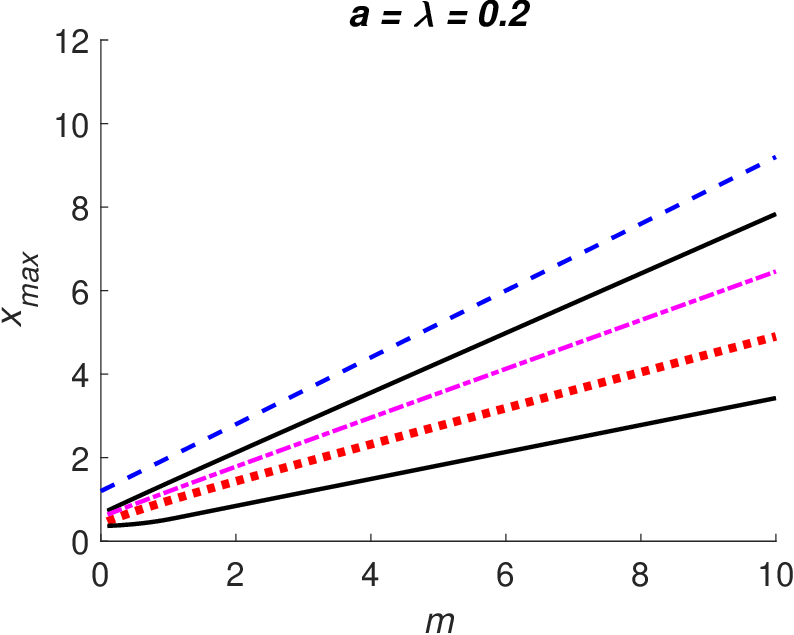}\hspace{0.5cm}
\includegraphics[scale = 0.5]{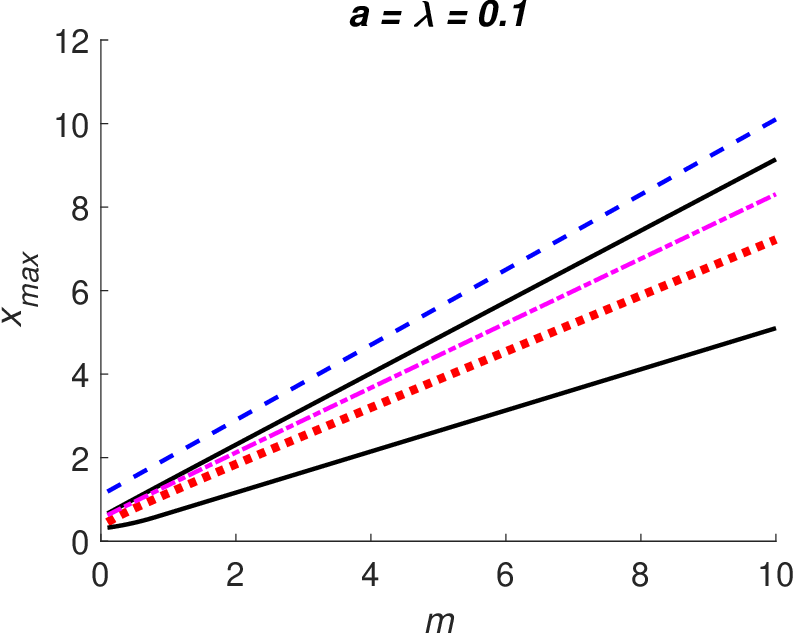}\vspace{1cm}
\includegraphics[scale = 0.5]{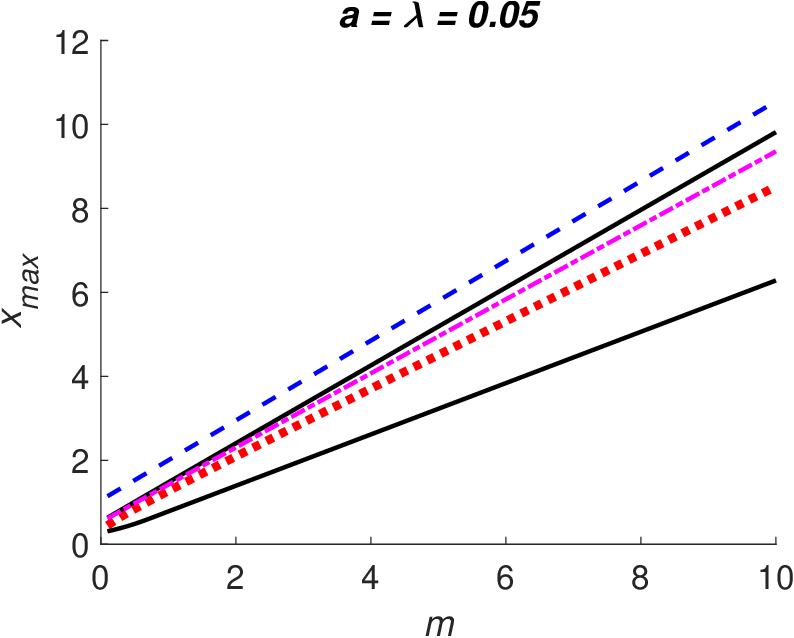}\hspace{0.5cm}
\includegraphics[scale = 0.5]{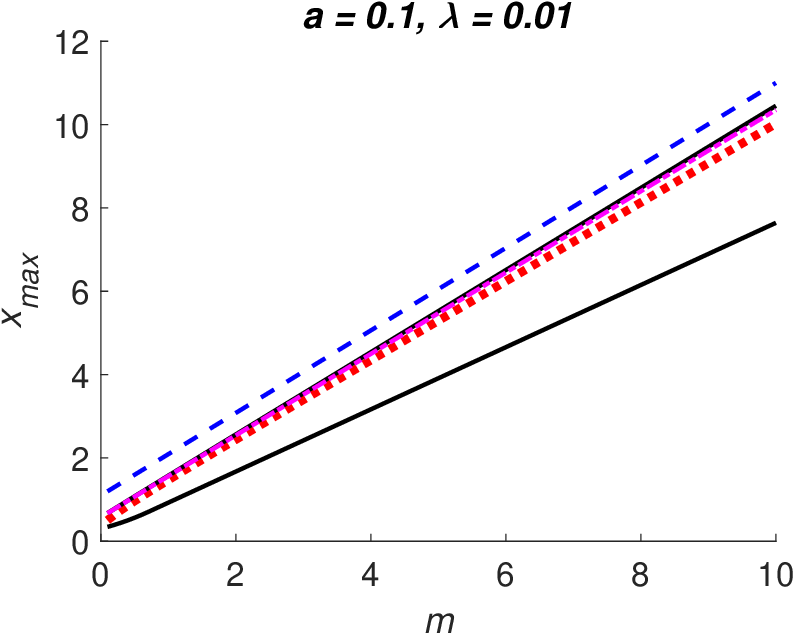}
\caption{
Estimates $\check x_1$ (with $s_0 = 0.8$) and $\hat x_1$ (black, solid), the estimate in Remark \ref{re:better} (magenta, dash-dot),
and the estimate in \eqref{eq:simple_x_max} (blue, dashed).
The red dotted curves give $x_{max}$ on (simulations of) the stable limit cycle. }
\label{fig:xmax}
\end{center}
\end{figure}

\subsection{Numerical results in Region 1}
\label{sec:numerics_reg1}

To achieve accurate numerics of system \eqref{eqmain}
under assumption \eqref{eq:assumption_parameters} we recommend transforming the equations (e.g. log transformations)
to avoid variables taking on very small values.
Imposing linear approximations near the unstable equilibria at
$(x,s) = (0,0)$ and $(x,s) = (0,1)$ are also helpful.
Indeed, using a standard ode-solver directly on system \eqref{eqmain} in cases of \eqref{eq:assumption_parameters}
may result in misleading trajectories,
showing a far to large minimum predator and prey biomass on the cycle,
unless tolerance settings are forced to minimum values.


We will compare the estimates for $x_1$ in Statement 1 with simulations of the stable limit cycle.
In addition, we make comparison to the upper estimate in Remark \ref{re:better} as well as to the simple linear estimate
\begin{align}\label{eq:simple_x_max}
x_{max} \leq 1 + a + m (1 - \lambda),
\end{align}
which is an immediate upper estimate of $\hat x_1$ in Statement 1.
(This bound is given by the escaping eigenvector of the saddle at (0,1), see the proof of Lemma \ref{th:uppis}.)
The behaviour of $x_{max}$ is approximately linear in $m$, as can be seen in
Figure~\ref{fig:xmax} showing $x_{max}$ on simulations of the stable limit cycle together with the above mentioned estimates as
functions of $m$ for some values of $a$ and $\lambda$.
We remark that we use $s_0 = 0.8$ in the lower estimates, which seems to be allowed also for the cases of $a, \lambda$ not satisfying \eqref{eq:assumption_parameters}, see Fig. \ref{fig:smax}.

\section{Estimates in Region 2 and Region 3}
\label{sec:lower-estimates}

\setcounter{theorem}{0}
\setcounter{equation}{0}

In order to state and prove our estimates for trajectories of system \eqref{eqmain} in Regions 2 and 3 we
define a trajectory $T^*$ as the solution curve with initial condition $(u,\lambda^*)$ with $u > h(\lambda^*)$ where $\lambda^* \leq \lambda$.
The trajectory $T^*$ crosses the isocline $x = h(s)$ at a minimal $s$-value at a point which we denote by $(h(s_u), s_u)$ and then it intersects $s = \lambda^*$ at a point which we denote by $(v,\lambda^*)$.
For these intersections we have the following estimates.

\begin{lemma}\label{th:minimal-lemma}
For the intersection of trajectory $T^*$ with the isocline $x = h(s)$ at $(h(s_u), s_u)$ it holds that
\begin{align}\label{eq:sminest}
\lambda^*\, e^{ - \frac{1}{m \lambda^*}\left(u - a - a \ln\frac u a\right) - 1} \,< \, s_u \, <\,  \lambda^* \, e^{- \frac{1}{m\lambda^*}\left( u - a - h(\lambda^*)\ln\frac{ u}{a}\right) },
\end{align}
and at the intersection with the isocline $s \,=\, \lambda^*$ at $(v,\lambda^*)$ it holds that
\begin{equation}\label{e3c}
\check z_1 \, u \, e^{ -\frac{u}{a}} \, < \, v \, <\,  \hat z_2 \, u\, e^{-\frac{u}{h(\lambda^*)}} \,< \, \hat z_0 \, u\, e^{-\frac{u}{h(\lambda^*)}}
\end{equation}
where $\check z_1 = z_1(u/a)$, $\hat z_2 = z_2(u/h(\lambda^*))$, $\hat z_0 = z_0(u/h(\lambda^*))$, and $1 < z_1, z_2, z_0 < e$ are decreasing functions defined in display \eqref{zi} below.
\end{lemma}

\noindent
{\it Proof of Lemma \ref{th:minimal-lemma}.}
We will trap trajectory $T^*$ between an upper and a lower barrier which we construct as follows.
Recall the phase-plane equation \eqref{eq:phase} and observe that
since $2\lambda + a < 1$ by assumption we have $a < h(s) < h(\lambda)$ in Region 2 and Region 3.
Therefore
\begin{align}\label{eq:tjosan}
\frac{a - x}{x (s - \lambda)}
>
\frac{m}{s}\frac{ds}{dx}
>
\frac{h(\lambda) - x}{x (s - \lambda)}
\end{align}
which produces separable ODEs.
Let
\begin{align*}
U_H(x,s) = \frac{1}{m} (x - H \ln x) + s - \lambda \ln s
\end{align*}
and observe that from \eqref{eq:tjosan} it follows that trajectory $T^*$ must stay below the level curve $\bar s(x)$ to $U_{h(\lambda^*)}$,
and above the level curve $\underbar s(x)$ to $U_a$,
whenever $s < \lambda^*$, where the barriers $\bar s(x)$ and $\underbar s(x)$ are defined through
\begin{align*}
U_{h(\lambda^*)}(x, \bar s(x)) = U_{h(\lambda^*)}\left(u,\lambda^*\right) \quad \text{and} \quad U_{a}(x, \bar s(x)) = U_{a}\left(u,\lambda^*\right).
\end{align*}

We now derive the estimates \eqref{eq:sminest} for the minimal prey biomass.
We begin with the upper estimate, for which we will find an upper estimate of the solution $\hat s$ to
\begin{align*}
U_{h(\lambda^*)}(h(\hat s), \hat s) = U_{h(\lambda^*)}(u, \lambda^*).
\end{align*}
%
The above equation is equivalent to
\begin{align*}
h(\hat s) - h(\lambda^*) \ln h(\hat s) + m(\hat s - \lambda^* \ln \hat s) = u - h(\lambda^*)\ln u  + m(\lambda^* -\lambda^* \ln \lambda^*).
\end{align*}
Since the function $x - A \ln x$ is decreasing for $0<x<A$, $\hat s < \lambda^*$ and $h(\hat s) \in (a, h(\lambda^*))$  we still get an upper estimate of $s_u$ by
$\hat s$ if we let it be the solution of
\begin{align*}
a - h(\lambda^*) \ln a + m(\hat s - \lambda^* \ln \hat s) = u - h(\lambda^*)\ln u + m(\lambda^* -\lambda^* \ln \lambda^*)
\end{align*}
which is equivalent to
\begin{align*}
 \hat s - \lambda^* \ln \hat s = \frac{1}{m}\left(u - h(\lambda^*)\ln u + m(\lambda^* -\lambda^* \ln \lambda^*) - a + h(\lambda^*) \ln a\right).
\end{align*}
We estimate from above once more by letting $\hat s$ solve
\begin{align*}
 \lambda^* - \lambda^* \ln \hat s = \frac{1}{m}\left(u - h(\lambda^*)\ln u + m(\lambda^* -\lambda^* \ln \lambda^*) - a + h(\lambda^*) \ln a\right)
\end{align*}
which gives
\begin{align*}
  \hat s = \lambda^* e^{- \frac{1}{m \lambda^*}\left(u - a - h(\lambda^*)\ln\frac u a\right) }
\end{align*}
and thus the upper estimate in \eqref{eq:sminest} is proven.

We proceed with the lower estimate,
for which we will find a lower estimate of the solution $\check s$ to
\begin{align*}
U_{a}(h(\check s), \check s) = U_{a}(u, \lambda^*)
\end{align*}
giving
\begin{align*}
h(\check s) - a \ln(h(\check s)) + m(\check s - \lambda^* \ln \check s) = u - a\ln u + m(\lambda^* -\lambda^* \ln \lambda^*).
\end{align*}
Since the function $x - A \ln x$ is decreasing for $0 < x < A$, increasing for $0 < A < x$,
$\check s < \lambda^*$ and $h(\check s) \in (a, h(\lambda^*))$ we still get a lower estimate of $s_u$ by
$\check s$ if we let it be the solution of
\begin{align*}
a - a \ln a + m(\check s - \lambda^* \ln \check s) = u - a \ln u + m(\lambda^* -\lambda^* \ln \lambda^*)
\end{align*}
which is equivalent to
\begin{align*}
 \check s - \lambda^* \ln \check s = \frac{1}{m}\left(u - a \ln u + m(\lambda^* -\lambda^* \ln \lambda^*) - a + a \ln a\right).
\end{align*}
We estimate from below once more by letting $\check s$ solve
\begin{align*}
 - \lambda^* \ln \check s = \frac{1}{m}\left(u - a \ln u + m(\lambda^* -\lambda^* \ln \lambda^*) - a + a \ln a\right)
\end{align*}
which gives
\begin{align*}
  \check s = \lambda^* e^{- \frac{1}{m \lambda^*}\left(u - a - a \ln\frac u a\right) - 1}
\end{align*}
and thus the lower estimate in \eqref{eq:sminest} is proven.
This completes the proof of \eqref{eq:sminest}.

\bigskip

We now prove estimate \eqref{e3c}.
For the lower bound we solve
\begin{align*}
U_{a}(\check x, \lambda^*) = U_{a}\left(u,\lambda^*\right)
\end{align*}
for $\check x$, which is equivalent to
\begin{align}\label{eq:checkx}
\check x - a \ln \check x = u - a \ln u.
\end{align}
Likewise, for the upper bound we solve $\hat x$ out of the equation
\begin{align*}
U_{h(\lambda^*)}(\hat x, \lambda^*) = U_{h(\lambda^*)}\left(u,\lambda^*\right),
\end{align*}
that is,
\begin{align}\label{eq:hatx}
\hat x - h(\lambda^*) \ln \hat x = u - h(\lambda^*) \ln u.
\end{align}
It then follows that $\check x < v < \hat x$ and we have estimate \eqref{e3c} if we can find the desired estimates of the solutions to
\eqref{eq:checkx} and \eqref{eq:hatx}.
To find such estimates we will use the following analytical approximations of Lotka-Volterra integrals from \cite{LS23},
to which we refer the reader for a proof.

\begin{lemma}\label{le:julelemma}
The solution $x < 1$ of the equation $x- \ln x = y - \ln y$ where $y > 1$ satisfies the relation $x=z\, Y$,
where $Y=y\, e^{-y}$, $1<z<e$ and $z = z(y)$ is decreasing in $y$.

Moreover,
let $d_i=e-1-c_i\, e$, $c_0=0$,
$c_1=\frac{e-2}{e-1}$, $c_2=\frac{1}{e}$, for $i=0,1,2$,
$Y=y\, e^{-y}$ and
\begin{equation}\label{zi}
z_i(y) =\frac{1-d_i\, Y-\sqrt{(1-d_i\, Y)^2 - 4 c_i\, Y}}{2c_i\, Y},\quad i = 1,2, \quad \text{and}\quad z_0(y)=\frac{1}{1-(e-1)Y}.
\end{equation}
Then the inequalities $1 < z_1 < z < z_2 < z_0 < e$ hold for $z$.
\end{lemma}

Observe that $\hat{x}$ is the solution to $\hat{x} - H \ln \hat{x} = u - H \ln u$, where $H = h(\lambda^*)$ and
by rescaling as $X = \hat{x}/H$ and $U = u/H$ we have $X - \ln X = U - \ln U$.
By assumption $u > h(\lambda^*)$ and thus $U > 1$. Therefore we can apply Lemma~\ref{le:julelemma} to obtain
$$
X = z(U)\, U\, e^{-U}
$$
where $z$ is a decreasing function of $U$.
As the same argument applies to the lower barrier $\check{x}$ we conclude
$$
z\left(\frac u a\right) u e^{-\frac{u}{a}} = \check{x} < v < \hat{x} = z\left(\frac{u}{h(\lambda^*)}\right) u e^{-\frac{u}{h(\lambda^*)}}.
$$
Finally, applying Lemma~\ref{le:julelemma} for estimating $z$ as $z_1 < z < z_2 < z_0$ we conclude
$$
z_1\left(\frac u a\right) u e^{-\frac{u}{a}} < v < z_2\left(\frac{u}{h(\lambda^*)}\right) u e^{-\frac{u}{h(\lambda^*)}} < z_0\left(\frac{u}{h(\lambda^*)}\right) u e^{-\frac{u}{h(\lambda^*)}}
$$
which is \eqref{e3c}. The proof of Lemma \ref{th:minimal-lemma} is complete. $\hfill\Box$\\

We next summarize implications for trajectory $T$ in the following two statements,
which are immediate consequences of Statement 1 and Lemma \ref{th:minimal-lemma} with $\lambda^* = \lambda$, for trajectory $T$.\\

\noindent
{\bf Statement 2}
\emph{
For the intersection of trajectory $T$ with isocline $x = h(s)$ at $P_2 = (h(s_2),s_2)$ it holds that
\begin{align*}
\lambda \,e^{- \frac{1}{m \lambda}\left(\hat x_1 - a - a\ln\frac{\hat x_1}{a}\right) - 1 } \,< \, s_2 \, < \, \lambda \, e^{- \frac{1}{m \lambda}\left(\check x_1 - a - h(\lambda)\ln\frac{\check x_1}{a}\right) }
\end{align*}
in which $\check x_1$ and $\hat x_1$ are as in Statement 1.}\\

\noindent
{\bf Statement 3}
\emph{
For the second intersection of trajectory $T$ with isocline $s = \lambda$ at $P_3 = (x_3,\lambda)$ it holds that
\begin{equation*}
 \check z_1 \, \hat x_1 \, e^{ -\frac{\hat x_1}{a}} \,<\, x_3 \,<\, \hat z_2 \, \check x_1 \, e^{-\frac{\check x_1}{h(\lambda)}} 
\end{equation*}
in which $\check x_1$ and $\hat x_1$ are as in Statement 1, $\check z_1 = z_1(\hat x_1/a)$, $\hat z_2 = z_2(\check x_1/h(\lambda))$ 
and
$1 < z_1, z_2 < e$ are the decreasing functions in display \eqref{zi}.}\\

\begin{figure}[h]
\begin{center}
\includegraphics[scale = 0.5]{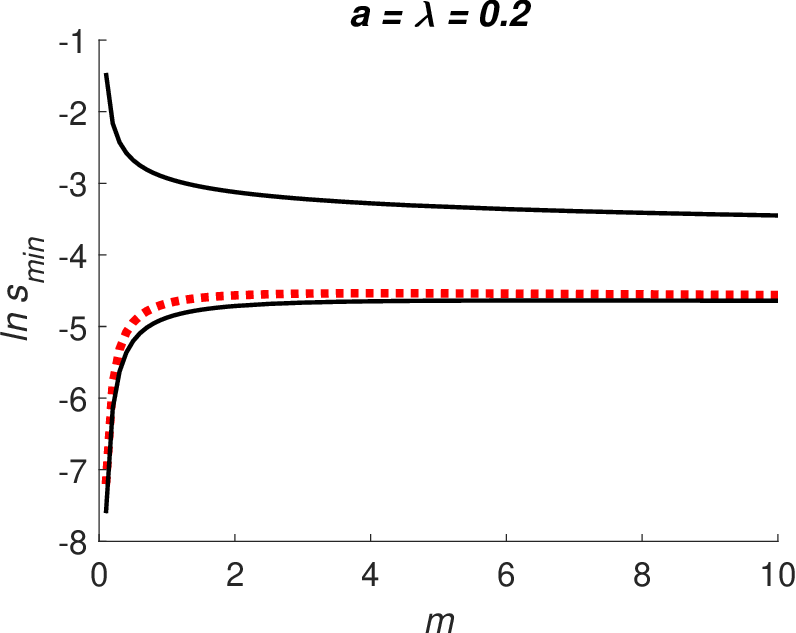}\hspace{0.5cm}
\includegraphics[scale = 0.5]{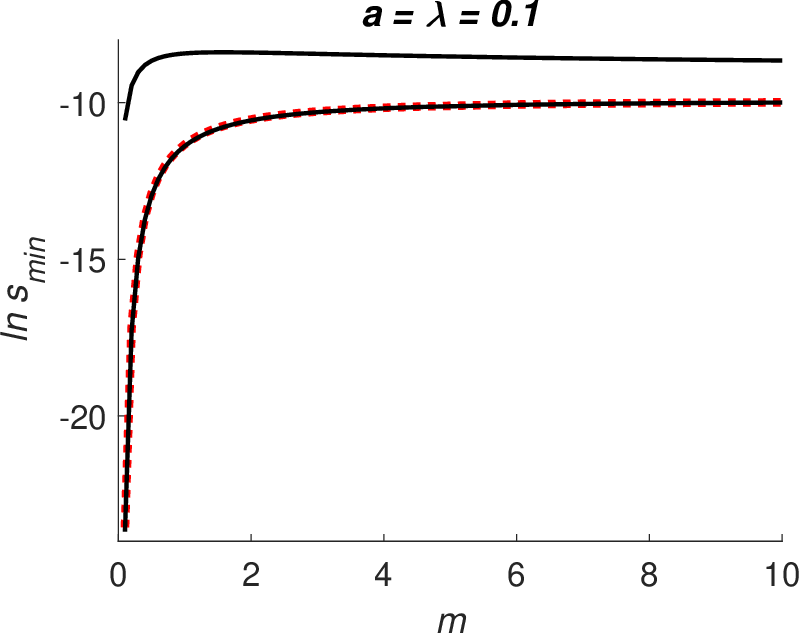}\vspace{1cm}
\includegraphics[scale = 0.5]{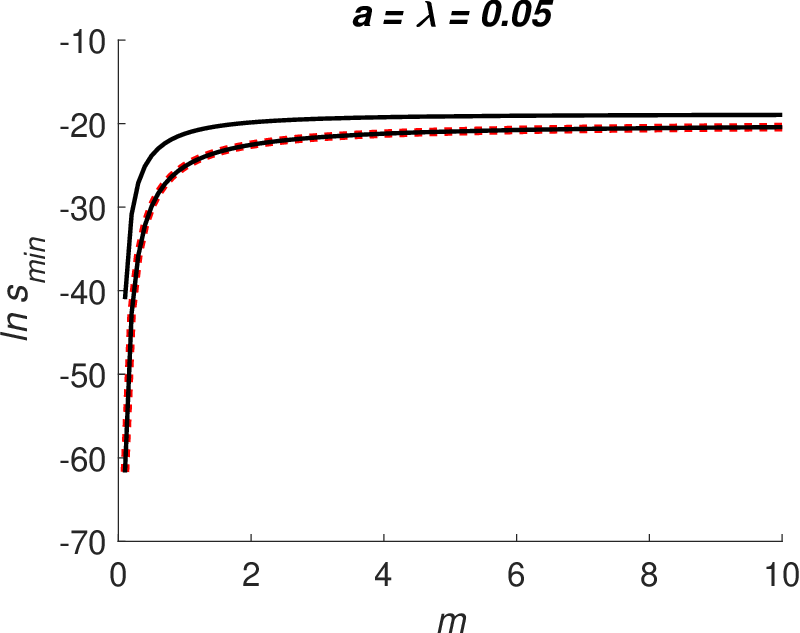}\hspace{0.5cm}
\includegraphics[scale = 0.5]{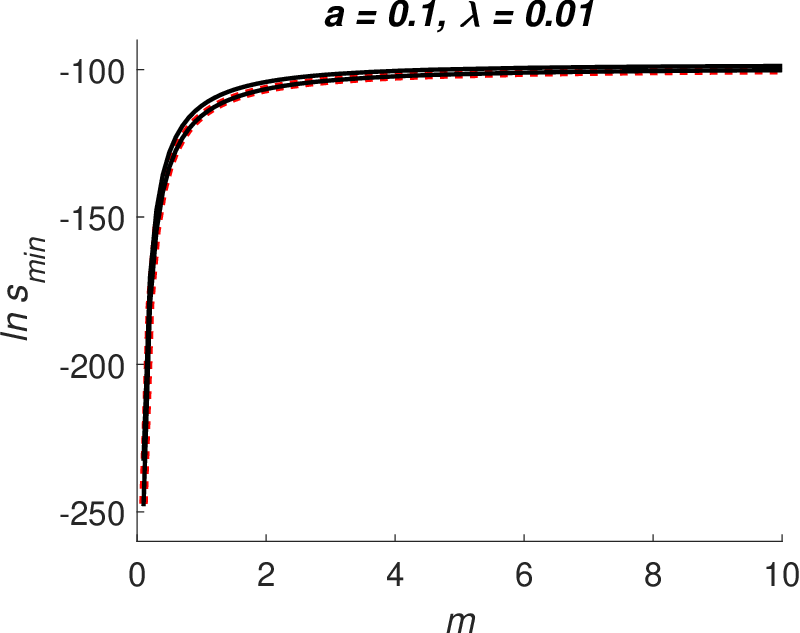}
\caption{
The estimates in Statement 2 (black, solid) together with $\ln s_{min}$ on (simulations of) the stable limit cycle.}
\label{fig:smin}
\end{center}
\end{figure}
\begin{figure}[h]
\begin{center}
\includegraphics[scale = 0.5]{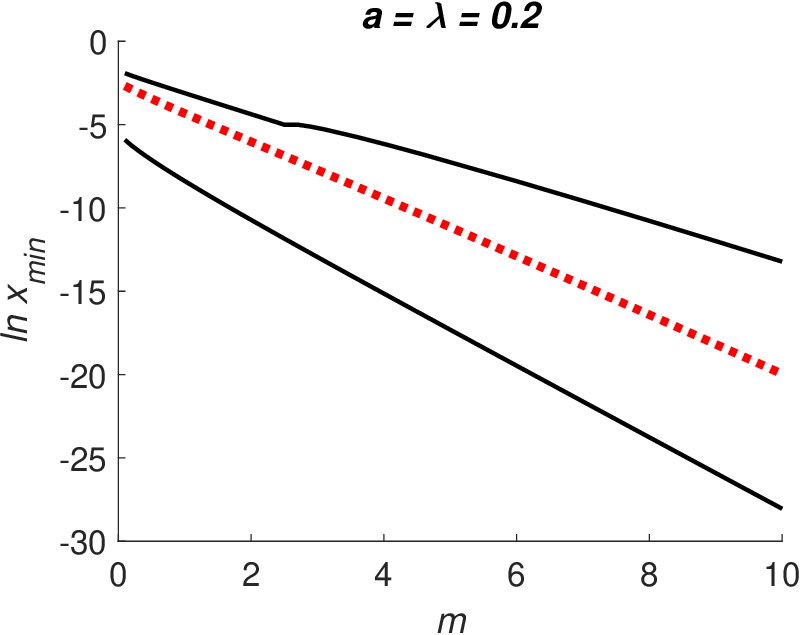}\hspace{0.5cm}
\includegraphics[scale = 0.5]{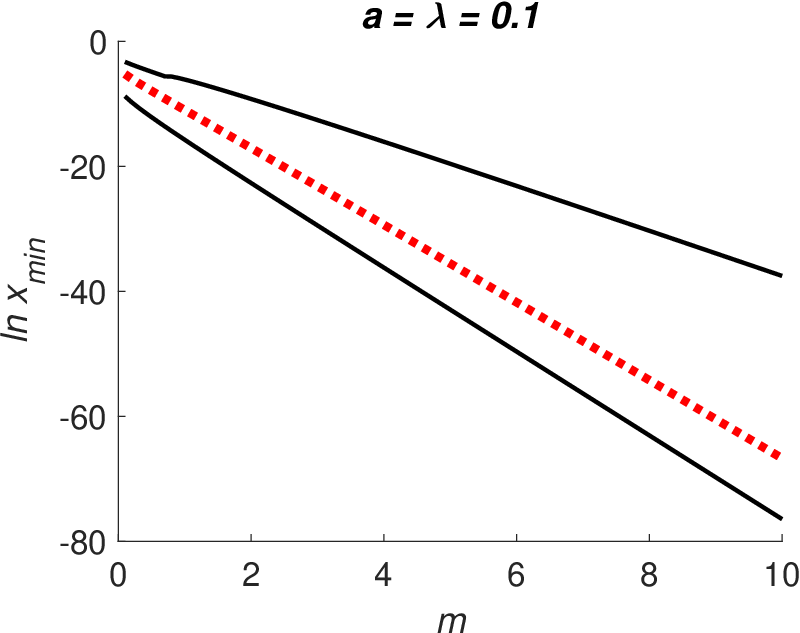}\vspace{1cm}
\includegraphics[scale = 0.5]{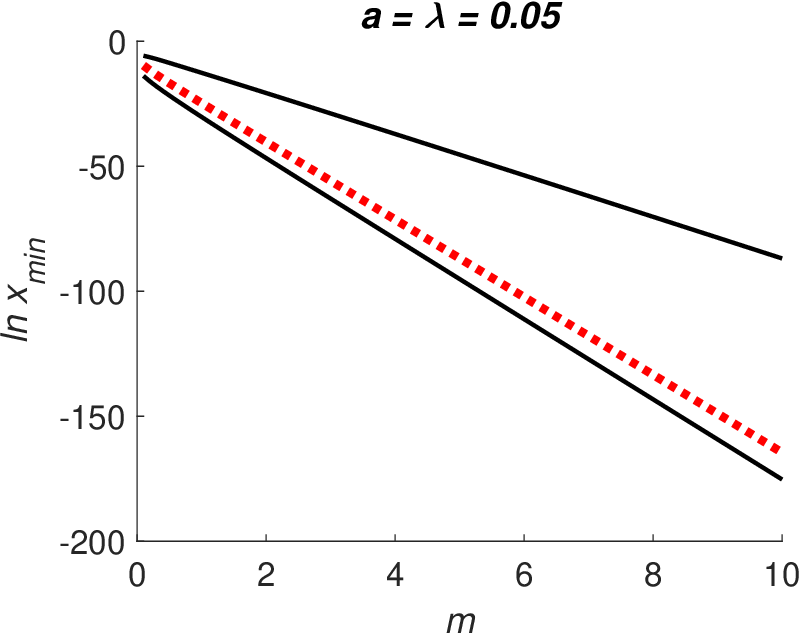}\hspace{0.5cm}
\includegraphics[scale = 0.5]{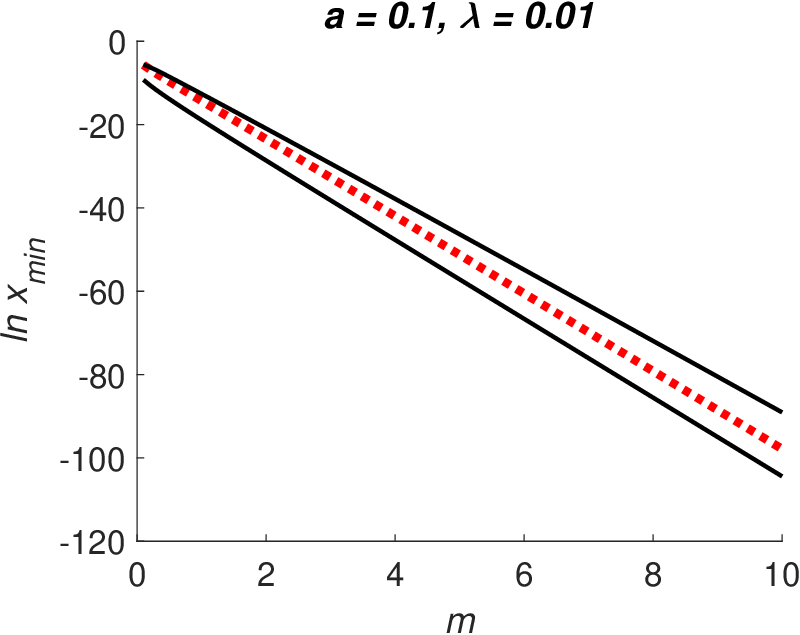}
\caption{
The estimates in Statement 3 (black, solid) together with $\ln x_{min}$ on (simulations of) the stable limit cycle.}
\label{fig:xmin}
\end{center}
\end{figure}

\subsection{Numerical results in Region 2 and Region 3}
\label{sec:numerics_reg23}

As in subsection \ref{sec:numerics_reg1} we end the section with numerical simulations.
Figure~\ref{fig:smin} (\ref{fig:xmin}) show $\ln s_{min} (\ln x_{min})$ on simulations of the stable limit cycle together with the estimates in Statement 2  (Statement 3) as
functions of $m$ for the same parameter values as in Figure \ref{fig:xmax}.
We can observe that the lower limit is good for small $\lambda$.



\section{Estimates in Region 4} 

\label{sec:reg4}

\setcounter{theorem}{0}
\setcounter{equation}{0}

In this section we prove the following statement for trajectory $T$:\\

\noindent
{\bf Statement 4}
\emph{Suppose that assumption \eqref{eq:assumption_parameters} holds and that $s_0 \geq 0.8$.
Then trajectory $T$, after intersecting $s\,=\,\lambda$ at $P_3$, intersects isocline $x = h(s)$ at $P_4 = (h(s_4),s_4)$ with $s_4 > 0.8$.}\\

\noindent
We split the proof into two steps of which the first considers estimates in the region $s \in (\lambda, 0.7]$
and the second considers estimates in $s \in (0.7, 1)$.
The first step, of which the proof is more lengthy and given at the end of the section,
consists in proving the following estimate for the $x$-value of the intersection of trajectory $T$ with $s = 0.7$:

\begin{lemma}\label{le:step1}
Suppose that assumption \eqref{eq:assumption_parameters} holds, that $s_0 \geq 0.8$ and $s_\gamma = 0.7$.
Then trajectory $T$, after intersecting $s\,=\,\lambda$ at $P_3$, intersects $s_\gamma$ at $(x_\gamma,s_\gamma)$ in which
\begin{align}\label{eq:etasimple}
x_\gamma \leq \hat \eta(m), \quad \text{for}\quad 0\leq m < \infty,
\end{align}
where $\hat \eta = \hat \eta (m)$ is given by \eqref{eq:etatjohej1} when $a,\lambda \leq 1/20$
and by \eqref{eq:etatjohej11} when $a \leq 1/10$ and $\lambda \leq a/10$.
\end{lemma}
The second step consists in proving that, given the estimate \eqref{eq:etasimple},
the trajectory $T$ intersects the isocline $x = h(s)$ at an $s$-value greater than 0.8. 
%
%
To prove this second step we establish the following lemma, which we state and prove
in a bit more general setting.

\begin{lemma}\label{le:step2}
Suppose that $s_\gamma \in (0,1)$, $M \in (0,s_\gamma]$, $x_\gamma \in (0,M(1-s_\gamma))$ and that $a, \lambda, m \in [0,\infty)$.
Then the trajectory of system \eqref{eqmain} starting at $(x_\gamma, s_\gamma)$ intersects isocline $x = h(s)$ at $(h(s_4),s_4)$ satisfying
$$
s_4 > 1 - \left( {\left(-d\right)^{m} x_\gamma^{M}} \;\frac{\left(m+M\right)^{m}}{M^{M} m^{m}} \right)^{\frac{1}{M+m}} \quad \text{where} \quad
d\,=\,s_\gamma+\frac{x_\gamma}{m + M} - 1.
$$
\end{lemma}

\noindent
{\it Proof}.
We consider trajectories of system \eqref{eqmain} in region
$$
E_{M,s_\gamma} \,=\,\lbrace \left(x,s\right) \vert \, 0\,<\,x\,<\,M\left(1-s\right), s\,>\,s_\gamma\rbrace ,
$$
where $M\,\leq\, s_\gamma$.
Observe that
$$
M\left(1-s\right)\,\leq\, s_\gamma\left(1-s\right)\,<\,s\left(1-s\right)\,<\,\left(s+a\right)\left(1-s\right)\,=\,h\left(s\right).
$$
In $E_{M,s_\gamma}$ we get the estimates
\begin{equation}
s'\,>\,\left(M\left(1-s\right)-x\right)s, \quad x'\,<\,m s x \quad \text{and hence} \quad \frac{ds}{dx} \,>\, \frac{M\left(1-s\right) - x}{m x}.
\label{5a}
\end{equation}
Let us consider a trajectory with initial condition $x\left(0\right)\,=\,x_\gamma,\, s\left(0\right)\,=\,s_\gamma$, where $x_\gamma\,<\,M\left(1-s_\gamma\right)$.
Using \eqref{5a},
we conclude that as long as this trajectory remains in $E_{M,s_\gamma}$,
it will be in the subregion bounded by the trajectory of the linear system
\begin{equation}
s'\,=\, M \left(1-s\right)-x,\quad x'\,=\,m x,
\label{5b}
\end{equation}
with initial condition  $x\left(0\right)\,=\,x_\gamma, s\left(0\right)\,=\,s_\gamma$ and the lines $x\,=\,M\, \left(1-s\right)$ and $x\,=\,0$.
Solving system \eqref{5b} we find that the trajectory follows the curve
\begin{equation}
s\,=\,d\left( \frac{x_\gamma}{x} \right)^\frac{M}{m} +1 -\frac{x}{m + M} \quad \text{with} \quad
d\,=\,s_\gamma+\frac{x_\gamma}{m + M} - 1.
\label{5c}
\end{equation}
The trajectory leaves  $E_{M,s_\gamma}$ when $x\,=\,M\, \left(1-s\right)$.
(Observe that then $s'\,=\,0$ for \eqref{5b}).
Substituting $x\,=\,M\, \left(1-s\right)$ into \eqref{5c} we get
$$
\frac{d x_\gamma^\frac{M}{m}}{M^\frac{M}{m}\left(1-s\right)^\frac{M}{m}} +1-s -\frac{M}{m+M} \left(1-s\right) \,=\, 0,
$$
which, if $d < 0$, is equivalent to
\begin{equation}\label{5d}
1-s \,=\, \left( {\left(-d\right)^{m} x_\gamma^{M}} \;\frac{\left(m+M\right)^{m}}{M^{M} m^{m}} \right)^{\frac{1}{M+m}}.
\end{equation}
The above expression for $1-s$ increases with $x_\gamma$,
for all $m > 0$ and $M \,>\,\frac{x_\gamma}{1-s_\gamma}$, because
\begin{align*}
\frac{\partial}{\partial x_\gamma} \left(\left(-d\right)^m x_\gamma^M\right) &\,=\, - x_\gamma^{M - 1} \left(-d\right)^{m - 1}\left( \frac{x_\gamma m}{m+M} + M\left(s_\gamma +\frac{x_\gamma}{m+M} -1\right) \right)\\ &\,=\, x_\gamma^{M - 1} \left(-d\right)^{m - 1}\left(M\left(1-s_\gamma\right)-x_\gamma\right) \,>\, 0.
\end{align*}
Thus, a lower boundary for the maximal $s$ of the trajectory of system \eqref{eqmain} starting at $(x_\gamma,s_\gamma)$ can be calculated from equation \eqref{5d}.
The proof of Lemma~\ref{le:step2} is complete. $\hfill\Box$\\

Before going into the proof of the first step, i.e. Lemma \ref{le:step1},
we show how Statement 4 follows from Lemma \ref{le:step1} and Lemma \ref{le:step2}.\\

\noindent
{\it Proof of Statement 4}.
We denote by $\alpha$ the upper estimate of $1-s_{max}$ with
\begin{equation*}
  \alpha = \alpha_1\, \alpha_2\,  \alpha_3,\quad
   \alpha_1 = \delta^{\frac{m}{m+M}},\quad
 \alpha_2 = \left(\frac{x_\gamma}{M}\right)^{\frac{M}{m+M}},\quad
  \alpha_3 = \left(\frac{m+M}{m}\right)^{\frac{m}{m+M}},
\end{equation*}
\noindent
where $\delta =1-s_\gamma-\frac{x_\gamma}{m+M}$.
We consider $\alpha$ for $m>0$ and take $M = s_\gamma = 0.7$.
We will make use of the following lemma:


\begin{lemma}\label{le:toppen}
Expressions $\alpha_1, \alpha_2$ and $\alpha_3$ satisfy the following properties:
\begin{itemize}
\item[A:] $\alpha_3$ has maximum $e^\frac{1}{e}$ at $m=\frac{M}{e-1}$,
is increasing for $m<\frac{M}{e-1}$ and decreasing for $m>\frac{M}{e-1}$.
%
%
%
\item[B:] $\alpha_2$ is increasing for $0<m<0.3$ and $0.3<m<m_1$ and decreasing for $m>m_1$,
where $m_1$ is defined below.
\item[C:] $\alpha_1$ is decreasing for $0<m<0.3$ and $m>0.3$.
\end{itemize}
\end{lemma}
%
%
Using Lemma \ref{le:toppen} for $m>1$ we notice that the value of  $\alpha_3$ is less
than $e^\frac{1}{e}$ and the value of  $\alpha_2$ is less than its value for $m=m_1$ and
$\alpha_1$ is less than its value for $m=1$. From this we easy calculate that
the product $\alpha <0.2$.

For $0.5<m<1$ we notice that the value of  $\alpha_3$ is less
than $e^\frac{1}{e}$ and the value of  $\alpha_2$ is less than its value for $m=1$ and
$\alpha_1$ is less than its value for $m=0.5$. From this we find $\alpha <0.2$.

For $0.3<m<0.5$ we notice that the value of  $\alpha_3$ is less
than $e^\frac{1}{e}$ and the value of  $\alpha_2$ is less than its value for $m=0.5$ and
$\alpha_1$ is less than its value for $m=0.3$. From this we find $\alpha <0.2$.

For $0.25<m<0.3$ we notice that the value of $\alpha_3$ is less
than its value for $m=0.3$ and the value of $\alpha_2$ is less than its value for $m=0.3$ and
$\alpha_1$ is less than its value for $m=0.25$. From this we find $\alpha <0.2$.

For $m<0.25$ we notice that the value of  $\alpha_3$ is less
than its value for $m=0.25$ and the value of $\alpha_2$ is less than its value for $m=0.25$ and
$\alpha_1$ is less than one. From this we find $\alpha <0.2$.

\bigskip

In the case when $a=0.1,\, \lambda=0.01$ the statement in Lemma \ref{le:toppen} still holds
with the only difference that $m_1\approx 3.06$ and the proofs are similar.
Thus, also in this case we can prove that $\alpha <0.2$ by examinations on interval.
We consider $\alpha_i,\, i=1.2,3$ as functions of $m$ and conclude:
\begin{align*}
\alpha <
\begin{cases}
    \alpha_1(1)\, \alpha_3(1)\, \max(\alpha_2)   < 0.193  \quad & \mbox{for} \quad 1 < m,\\
    \alpha_1(0.7)\, \alpha_2(1)\, \max(\alpha_3) < 0.194  \quad & \mbox{for} \quad 0.7 < m \leq 1,\\
    \alpha_1(0.5)\, \alpha_2(0.7)\, \max(\alpha_3) <0.19  \quad & \mbox{for} \quad 0.5 < m \leq 0.7,\\
    \alpha_1(0.3)\, \alpha_2(0.5)\, \max(\alpha_3) <0.187 \quad & \mbox{for} \quad 0.3 < m \leq 0.5,\\
    \alpha_1(0.28)\, \alpha_2(0.3)\, \alpha_3(0.3) <0.198 \quad & \mbox{for} \quad 0.28 < m \leq 0.3,\\
    \alpha_1(0.2)\, \alpha_2(0.28)\, \alpha_3(0.28) <0.2  \quad & \mbox{for} \quad 0.2 < m \leq 0.28,\\
    \alpha_1(0.1)\, \alpha_2(0.2)\, \alpha_3(0.2) <0.17   \quad & \mbox{for} \quad 0.1 < m \leq 0.2,\\
    1\cdot\, \alpha_2(0.1)\, \alpha_3(0.1) < 0.12         \quad & \mbox{for} \quad m \leq 0.1.
\end{cases}
\end{align*}









\noindent
From this we conclude that $\alpha <0.2$ and the proof of Statement 4 is complete. $\hfill\Box$

\bigskip

\noindent
\emph{Proof of Lemma \ref{le:toppen}.} Statement A is trivial.
To prove statement B we denote $u= \frac{\partial (\ln\alpha_2)}
{\partial m} \frac{(m+M)^2}{M}$. In this notation we get
$u= (M+m)\frac{x'_\gamma}{x_\gamma} - \ln \frac{x_\gamma}{M}$, where differentiation is taken with resp to $m$.
If $x_\gamma=(\bar a\, m +b)\, e^{-c\, m}$ we get
\begin{equation*}
u = (m+M)\frac{\bar a - b\, c - \bar a\, c\, m}{\bar a\, m+b} +c\, m -\ln (\bar a\, m +b) +\ln M.
\end{equation*}
Calculations give $u'=-\frac{\bar a^2\, (m+M)}{(\bar a\, m +b)^2}$
for the derivative with respect to $m$ and thus $u$ is decreasing in $m$.

For the actual values of $\bar a,\, b,\, c$, in the case $a=\lambda=0.05$,
$u$  has a root at $m = m_1\approx 4.11$ when $m > 0.3$.
For $m<0.3$, $x_\gamma<1$ is increasing in $m$ and $\frac{M}{m+M}$ is decreasing,
implying that $\alpha_2$  is increasing.

\bigskip

The derivative of $\ln\alpha_1$ with respect to $m$ is calculated as


\begin{equation*}
\frac{1}{\delta\, (m+M)^2} \left( M\,\delta\, ln\delta + m\, (m+M)\, \delta'\right)
\end{equation*}
\noindent
where
\begin{equation}
\delta' = - \left(\frac{x_\gamma}{m+M}\right)' =\frac{x_\gamma}{(m+M)^2} - \frac{x'_\gamma}{m+M} = -\frac{Q(m)\, e^{-c\, m}}
{(m+M)^2},
\label{deltadiff}
\end{equation}
\begin{equation*}
 Q(m)=-\bar a\,c\, m^2 - (b+\bar a\, M)\, c\, m + M\, (\bar a - b\, c) -b,
\end{equation*}
\noindent
and $x_\gamma=(\bar a\, m + b)\, e^{-c\, m}$.
For the actual values of $\bar a,\, b,\, c$ in the case $a=\lambda=0.05$ we get
$Q(m)>0$ for $m<0.3$ and thus from (\ref{deltadiff}) follows that $\delta'$ is
negative and $\alpha_1$ decreasing.

\bigskip

For $m>0.3$ we notice that $m\, (m+M)\, \delta' < M \delta\, | \ln\delta |$  and thus $\alpha_1$ is decreasing.
Indeed, from $Q(m)<Q(0.3)<0$ follows that $\delta'$ is positive and thus
$\delta\, \ln\delta$ is less than its value at $m=0.3$ (since $\delta < 1 - s_\gamma = 0.3 < e^{-1}$) which is less than -0.35 M.
For $m>1$ we also have $\frac{m}{m+M}<1$ and $m\, (m+M)\, \delta' < m^2\, e^{-2\, m}<e^{-2}<0.14 $
and for $0.3<m<1$ we get $\frac{m}{m+M}<1$ and $m\, (m+M)\, \delta' < m\, e^{-2\, m}<\frac1{2e}<0.2 $.
Recalling $M = 0.7$ the proofs are finished.$\hfill\Box$

\bigskip


We now turn to the first step (Lemma \ref{le:step1}) which proof is based on a number of auxiliary lemmas.
In order to state the first of these lemmas we introduce,
for parameters $k, s_\gamma \in (0,1)$,
the function
\begin{equation}\label{eq:etadef}
\eta \left(a,\lambda,m\right) \,=\, \left( e^{\frac{\lambda}{s_\gamma}}\, \frac{s_\gamma+a}{1-s_\gamma}\, \frac{1}{a+\lambda}\right)^{\frac{m}{k}}\, x_3.
\end{equation}
%
Our first auxiliary lemma yields as follows.
\begin{lemma}\label{le:lemma18}
Suppose $k \in (0,1)$ and $s_\gamma \in (0.5,1)$. If
$
\eta \left(a,\lambda,m \right)\,<\,\left(1-k\right)\, h\left(s_\gamma\right)
$
and $x_3\,<\,\left(1-k\right)\, h\left(\lambda\right)$, then the trajectory $T$ intersects $s \,=\, s_\gamma$ before escaping Region~4 and at the intersection $\left(x_\gamma,s_\gamma\right)$ the estimate $x_\gamma\,<\,\eta \left(a,\lambda,m \right)$ holds.
\end{lemma}
We build the proof of Lemma \ref{le:lemma18} on the following two lemmas, which proofs we postpone until the end of the section.

\begin{lemma}\label{le:lemma20}
Suppose that $k \in (0,1)$ and $x_3\,<\,\left(1-k\right)\, h\left(\lambda\right)$.
Then as long as the trajectory with initial condition $(x_3,\lambda)$ stays in the region determined by
$x\,<\,\left(1-k\right)\, h\left(s\right)$ it holds that
\begin{equation}
x\,<\,x_3\, B^{\frac{m}{k}}
\label{e7a}
\end{equation}
where
$$
B = \left( \frac{s+a}{s}\right)^{\frac{\lambda}{a}}
\left( \frac{\lambda}{\lambda +a}\right)^{\frac{\lambda}{a}}
\left( \frac{\left(1-\lambda\right)\left(s+a\right)}{1-s}\right)^{\frac{1-\lambda}{1+a}}
\left( \frac{1}{\lambda +a}\right)^{\frac{1-\lambda}{1+a}}.
$$
Moreover, if $s \,\geq\, 0.5$ then $B \leq e^{\frac{\lambda}{s}} \, \frac{s+a}{1-s} \, \frac{1}{a+\lambda}$ and hence also
\begin{equation}
x\,<\,x_3\, B^{\frac{m}{k}}\,<\,x_3\, \left( e^{\frac{\lambda}{s}} \, \frac{s+a}{1-s} \, \frac{1}{a+\lambda} \right)^{\frac{m}{k}} = \eta(a,\lambda,m).
\label{e7b}
\end{equation}
\end{lemma}

\begin{lemma}\label{le:lemma21}
Suppose that $k \in (0,1)$, $s_\gamma \in (0.5,1)$, $x_3\,<\,\left(1-k\right)\, h\left(\lambda\right)$ and $x_3\, \left( B_\gamma \right)^{\frac{m}{k}}$
$\,<\, \left(1-k\right)\, h\left(s_\gamma\right)$,
where $B_\gamma$ is the value $B$ takes for $s\,=\,s_\gamma$. Then the trajectory with initial condition $(x_3,\lambda)$ intersects $x\,=\,\left(1-k\right)\, h\left(s\right)$ next time for $s\,>\,s_\gamma$ and stays inside the region $x\,<\,\left(1-k\right)\, h\left(s\right)$ before it intersects $s\,=\,s_\gamma$.
\end{lemma}

\noindent
{\it Proof of Lemma \ref{le:lemma18}}.
Since $s > 0.5$ we can conclude, from the estimate of $B$ derived in the proof of Lemma \ref{le:lemma20}, that
$
x_3\, \left( B_\gamma \right)^{\frac{m}{k}}\,<\, 
\eta(a,\lambda,m)
$.
%
Therefore, the assumptions in Lemma \ref{le:lemma18} imply the assumptions in Lemma \ref{le:lemma21} and it follows that
the trajectory $T$ will be inside the region $x\,<\,\left(1-k\right)\, h\left(s\right)$ until it intersects $s = s_\gamma$.
Thus, we can apply Lemma~\ref{le:lemma20} to obtain
\begin{align*}
x_\gamma < x_3\, \left( B_\gamma \right)^{\frac{m}{k}}\,<\,  \eta(a,\lambda,m),
\end{align*}
which is the sought after estimate.  $\hfill\Box$\\

Our next step is to build an upper estimate of the function $\eta(a,\lambda,m)$ defined in \eqref{eq:etadef}.
A first step is to apply our preceding estimates in Statement 3 for $x_3$ together with a lower estimate for $\check x_1$ to obtain
\begin{equation}\label{eq:etasimp1}
\eta \left(a,\lambda,m\right) \,\leq\, \left( e^{\frac{\lambda}{s_\gamma}}\, \frac{s_\gamma+a}{1-s_\gamma}\, \frac{1}{a+\lambda}\right)^{\frac{m}{k}}\, \hat z_2 \, \tilde x_1 \, e^{-\frac{\tilde x_1}{h(\lambda)}} := \bar \eta(a,\lambda,m).
\end{equation}
Here, $\hat z_2 = z_2(\tilde x_1/h(\lambda))$ where $1 < z_2 < e$ is the decreasing function in display \eqref{zi}
and $\tilde x_1 = c_0 + m c$ is a lower estimate of $\check{x}_1$ in Statement 1 in which
%
%
%
%
\begin{align}\label{eq:x1simp1}
c_0  &=
    \begin{cases}
        \frac{1}{4}     & \mbox{if $m < 0.3$}\\\notag
        h(0.8)          & \mbox{if $m \geq 0.3$}	
    \end{cases}\\
%
c &=
    \begin{cases}
        \frac{1-a_{max}}{2} -\lambda\, (1-\ln \lambda +\ln \frac{1-a_{max}}{2} )  & \mbox{if $m < 0.3$}\\
        0.8 -\lambda\, (1-\ln \lambda +\ln 0.8 )                     & \mbox{if $m \geq 0.3$}	
    \end{cases}
\end{align}
where $a_{max}$ is the highest allowed value of $a$.
To proceed with constructing an upper estimate of $\eta$ we need the following lemma,
which we prove in the end of the section.


\begin{lemma}\label{le:lemma19}
Let $\bar \eta = \bar \eta \left(a,\lambda,m\right)$ be as in \eqref{eq:etasimp1} but allowing for $z_i$, $i = 0,2$, (in place of $\hat z_2$)
where $z_0$ and $z_2$ is as defined in \eqref{zi} with $y = \frac{\tilde x_1}{h(\lambda)}$.
Then the derivatives of $\bar \eta$ with respect to $a$ and $\lambda$ are nonnegative
if $k \geq 1/3$, $a+\lambda \leq 0.2$, and if there exists $\kappa \in (0,1)$ so that $\kappa \,\tilde x_1 \geq a + \lambda$ and $(1-\kappa) k \ln \frac{s_0}{\lambda} \geq 1$.
\end{lemma}

\noindent
We are now ready to prove Lemma \ref{le:step1}. \\

\noindent
{\it Proof of Lemma \ref{le:step1}}.
We intend to establish \eqref{eq:etasimple}, i.e.
\begin{align*}
x_{\gamma} \leq \bar \eta \left(a,\lambda,m\right) \leq \hat \eta(m)
\end{align*}
for a suitable function $\hat \eta = \hat \eta(m)$.
We will begin with the second inequality, i.e. finding $\hat \eta$ as an upper estimate for the $\bar \eta$ function defined in \eqref{eq:etasimp1} and \eqref{eq:x1simp1}, and end by using Lemma \ref{le:lemma18} to conclude $x_\gamma  < \eta(a,\lambda,m)$.
We present the arguments in detail for the case $a, \lambda \leq 0.05$, and state only the explicit estimates in the case $a \leq 0.1, \lambda < 0.01$.

\bigskip

\noindent
{\bf The case $a,\lambda < 0.05$:}
Pick $k = 0.75$, $\kappa = 2/5$ and observe that \eqref{eq:x1simp1} implies $\tilde x_1 \geq 1/4$.
It is then easy to verify that all assumptions in Lemma \ref{le:lemma19} hold as long as assumption \eqref{eq:assumption_parameters} holds.
Indeed, $a+\lambda \leq 0.1$, $\kappa \tilde x_1 \geq 0.1$, and as $s_0 \geq 0.475$ we also see that $(1-\kappa) k \ln \frac{s_0}{\lambda} \geq 1$.
An application of Lemma \ref{le:lemma19} gives
\begin{align}\label{eq:estetatjohej}
\bar \eta \left(a,\lambda,m\right) \leq \bar \eta(0.05,0.05,m) \leq \hat{c}^{\frac{m}{k}} \hat z_2 \tilde x_1 e^{-\frac{\tilde x_1}{0.095}}
\end{align}
in which $\hat{c}$ is an upper roundoff of the value $e^{\frac{\lambda}{s_\gamma}}\, \frac{s_\gamma+a}{1-s_\gamma}\, \frac{1}{a+\lambda}$
takes with $s_\gamma = 0.7$ and $a = \lambda = 0.05$, and $\hat z_2$ and $\tilde x_1$ are evaluated at $a = \lambda = 0.05$.
To simplify further we observe that by \eqref{eq:x1simp1} and since $y = \frac{\tilde x_1}{h(\lambda)}$ we have $\frac{\partial y}{\partial m} > 0.$
Moreover, since $y > 1$ we also have
$
\frac{\partial Y}{\partial y} = (1-y) e^{-y} < 0
$
and, from \eqref{eq:deriv:z0z2>0} in the proof of Lemma \ref{le:lemma19} it follows that
\begin{align}\label{eq:monotonzinm}
\frac{\partial z_i}{\partial m} = \frac{\partial z_i}{\partial Y} \frac{d Y}{d y} \frac{\partial y}{\partial m} < 0.
\end{align}
It is thus possible to find an upper bound of $\hat z_i$, $i = 0,2$, valid in a range $m_0 \leq m$, by fixing $m = m_0$ in the $\hat z_i$s.
Therefore, let $c_z$ be the value $\hat z_2$ takes for $m = m_0$.
Then from \eqref{eq:estetatjohej} and \eqref{eq:monotonzinm} we have
\begin{align}\label{eq:etaest2}
\bar \eta \left(a,\lambda,m\right)
\leq \hat{c}^{\frac{m}{k}} c_z \tilde x_1 e^{-\frac{\tilde x_1}{0.095}}.
\end{align}
Moreover,
$
\frac{\partial}{\partial \tilde x_1} \left(\tilde x_1 e^{ - \frac{\tilde x_1}{0.095}} \right) < 0 
$
and hence we can approximate by replacing $\tilde x_1$ with $c_0 + m c$ where $c_0$ and $c$ are appropriate lower roundoffs for the expressions in \eqref{eq:x1simp1}. 
Continuing from \eqref{eq:etaest2} we end up with
\begin{align*}
\bar \eta \left(a,\lambda,m\right) \leq \hat{c}^{\frac{m}{k}} c_z \left( c_0 + m c \right) e^{-\frac{ c_0 + m c }{0.095}}
= c_z(c_0 + m c) e^{-\frac{c_0}{0.095} + m \left(\frac{\ln \hat c}{k} - \frac{c}{0.095}\right)}.
\end{align*}

When $m \geq m_0 = 0.3$ we obtain, with $k=0.75$: $c_z c_0 \leq 0.190$, $c_z c \leq 0.681$, $\frac{\ln \hat c}{k} - \frac{c}{0.095} \leq -2.048$, $\frac{c_0}{0.095} \geq 1.789$ and for $m \in (0, 0.3]$ we use $m_0 = 0$ and end up with $c_z c_0 \leq 0.324$, $c_z c \leq 0.404$, $\frac{\ln \hat c}{k} - \frac{c}{a+\lambda} \leq 1.099$, $-\frac{c_0}{a+\lambda} \leq -2.631$.
We conclude $\bar \eta(a,\lambda,m) \leq \bar \eta(0.05, 0.05, m) \leq \hat \eta(m)$ in which
\begin{align}\label{eq:etatjohej1}
\hat \eta(m) =
\begin{cases}
(0.324 + 0.404 m) e^{1.099 m - 2.631} & \quad \text{if} \quad  0 \leq m \leq 0.3\\
(0.190 + 0.681 m) e^{-2.048 m - 1.789} & \quad \text{if} \quad 0.3 < m.
\end{cases}
\end{align}
%


\bigskip

To complete the proof of Lemma \ref{le:step1} we will use Lemma \ref{le:lemma18} to ensure
$x_\gamma \leq \eta(a,\lambda,m)$ whenever \eqref{eq:assumption_parameters} holds.
To do so we have to ensure the validity of the assumptions in Lemma \ref{le:lemma18}.
To verify that $\eta(a,\lambda, m) \leq (1-k)h(s_\gamma)$ when $k = 0.75$ and $s_\gamma = 0.7$ it suffices to observe that
\begin{align*}
\hat \eta(m) \leq 0.05 \leq (1-k)h(s_\gamma), \quad \text{whenever} \quad m \in [0, \infty),
\end{align*}
$a, \lambda \leq 0.05$ and where $\hat \eta$ is from \eqref{eq:etatjohej1}.

To see that $x_3 \leq (1-k) h(\lambda)$ we first apply Statement 3 to obtain
\begin{align}\label{eq:firstxminestimatetjohej}
x_3 \leq \hat z_2 \tilde x_1 e^{- \frac{\tilde x_1}{h(\lambda)}}
\end{align}
and we let again $c_z$ be the value of $\hat z_2$ when $a = \lambda = 0.05$ and $m = m_0$.
Then $\hat z_2 \leq c_z$ whenever $m \geq m_0$ and \eqref{eq:assumption_parameters} hold because $\frac{\partial z_i}{\partial \lambda} > 0$ by \eqref{eq:deriv:2},
$\frac{\partial z_i}{\partial a} > 0$ by \eqref{eq:derivya} and \eqref{eq:deriv:2a},
and $\frac{\partial z_i}{\partial m} < 0$ by \eqref{eq:monotonzinm}.
Further, we let $c_x$ be the value of $\tilde x_1$ when $a = 0$, $\lambda = 0.05$ and $m = m_0$.
Then $c_x \leq \tilde x_1$ holds by monotinicity (recall \eqref{eq:x1simp1}),
and as $\tilde x_1 e^{-\frac{\tilde x_1}{h(\lambda)}}$ is decreasing in $\tilde x_1$ it follows from \eqref{eq:firstxminestimatetjohej} that
$$
x_3  \leq c_z c_x e^{- \frac{c_x}{h(\lambda)}}.
$$
For $m \geq m_0 = 0.3$ we obtain $c_z \leq 1.114$ and $c_x \geq 0.343$,
and for $m_0 = 0, m \in (0, 0.3)$ we have $c_z \leq 1.293$ and $c_x \geq 0.25$.
Therefore
\begin{align}\label{eq:secondxminestimatetjohej}
x_3 \leq
\begin{cases}
0.324\, e^{- \frac{1}{4 h(\lambda)}} & \quad \text{if} \quad  0 \leq m \leq 0.3\\
0.383\, e^{- \frac{0.343}{h(\lambda)}} & \quad \text{if} \quad 0.3 < m.
\end{cases}
\end{align}
For the right hand side we have
$
(1-k)h(\lambda) \geq 0.25h(\lambda)
$
as $k = 0.75$, which must exceed \eqref{eq:secondxminestimatetjohej} whenever $h(\lambda) \in (0, 0.095)$.
Plugging in $h(\lambda) = 0.095$ in both sides proves that it is true.
Thus, we conclude that
$x_3 \leq (1-k) h(\lambda)$.

\bigskip

\noindent
{\bf The case $a < 0.1, \lambda < 0.1 a$:}
We proceed as in the former case but now with $k = 2/3$, $\kappa = 1/2$, $a \leq 0.1$ and $\lambda \leq 0.01$.
We obtain $\bar \eta(a,\lambda,m) \leq \bar \eta(0.1, 0.01, m) \leq \hat \eta(m)$ in which
\begin{align}\label{eq:etatjohej11}
\hat \eta(m) =
\begin{cases}
(0.350 + 0.563 m) e^{1.113 m - 2.295} & \quad \text{if} \quad  0 \leq m \leq 0.3\\
(0.201 + 0.832 m) e^{-2.048 m - 1.652} & \quad \text{if} \quad 0.3 < m.
\end{cases}
\end{align}
It is easy to verify that $\eta(a,\lambda, m) \leq 0.08 \leq (1-k)h(s_\gamma)$.
Moreover,
\begin{align*}
x_3 \leq
\begin{cases}
0.350\, e^{- \frac{1}{4 h(\lambda)}} & \quad \text{if} \quad  0 \leq m \leq 0.3\\
0.428\, e^{- \frac{0.383}{h(\lambda)}} & \quad \text{if} \quad 0.3 < m,
\end{cases}
\end{align*}
and since
$
x_3 \leq \frac13h(\lambda)
$
whenever $h(\lambda) \in (0, 0.1089)$ we have also ensured $x_3 \leq (1-k) h(\lambda)$.
The proof of Lemma \ref{le:step1} is complete.
$\hfill\Box$\\

We now give proofs of the auxiliary Lemmas \ref{le:lemma20}, \ref{le:lemma21} and \ref{le:lemma19}.\\

\noindent
{\it Proof of Lemma~\ref{le:lemma20}}.
When $x\,<\,\left(1-k\right)\, h\left(s\right)$ we have $s'\,>\,k h\left(s\right) s\,>\,0$ and as $s > \lambda$ it also holds that $x'\,>\,0$.
Thus we get the inequality
\begin{equation*}
\frac{dx}{ds} \,<\, \frac{m x \left(s-\lambda\right)}{k h\left(s\right) s}.
\end{equation*}
Integrating gives
\begin{equation}\label{eq:F_def}
x\,<\, x_3 \left( \frac{F\left(s\right)}{F\left(\lambda\right)}\right)^{\frac{m}{k}}\quad \text{where}\quad
F\left(y\right)\,=\,\frac{\left(y+a\right)^{k_1}}{y^{k_2} \left(1-y\right)^{k_3}}
\end{equation}
and where $k_1\,=\,\frac{a+\lambda}{a\left(a+1\right)}$, $k_2\,=\,\frac{\lambda}{a}$ and $k_3\,=\,\frac{1-\lambda}{1+a}$.
But
\begin{equation*}
\frac{F\left(s\right)}{F\left(\lambda\right)}\,=\,
\left( \frac{s+a}{\lambda +a}\right)^{k_1}
\left( \frac{\lambda}{s}\right)^{k_2}
\left( \frac{1-\lambda}{1-s}\right)^{k_3}
\end{equation*}
and because $k_1\,=\,k_2+k_3$ we get
\begin{equation}\label{eq:EEEE}
\frac{F\left(s\right)}{F\left(\lambda\right)}\,=\,
\left( \frac{s+a}{s}\right)^{k_2}
\left( \frac{\lambda}{\lambda +a}\right)^{k_2}
\left( \frac{\left(1-\lambda\right)\left(s+a\right)}{1-s}\right)^{k_3}
\left( \frac{1}{\lambda +a}\right)^{k_3} = B,
\end{equation}
which, together with \eqref{eq:F_def}, proves estimate \eqref{e7a}.

To prove \eqref{e7b} we denote $B \,=\,  E_1 E_2 E_3 E_4$,
where $E_i, \, i\,=\,1,2,3,4$ are the factors in \eqref{eq:EEEE} defined in the obvious way.
We conclude
$$
E_1\,=\, \left(\left( 1+ \frac{1}{\frac{s}{a}}\right) ^\frac{s}{a} \right)^\frac{\lambda}{s} \,<\, e^\frac{\lambda}{s},\quad
E_2\,=\,\frac{1}{\left( 1+\frac{1}{k_2}\right)^{k_2}}\,<\,1,
$$
and by using $k_3 < 1$ and $s \geq 0.5$ for estimating $E_3$,
and $k_3 < 1$ and $a+\lambda < 1$ for estimating $E_4$, we get
$$
E_3\,<\, \left(\frac{s+a}{1-s}\right)^{k_3} \,<\, \frac{s+a}{1-s},\quad
E_4\,<\, \frac{1}{\lambda +a}
$$
%
All these estimates together with \eqref{eq:EEEE} give \eqref{e7b}. $\hfill\Box$\\

\noindent
{\it Proof of Lemma~\ref{le:lemma21}}.
We first claim that the assumptions in the lemma imply
\begin{equation}\label{eq:lemma21_claim}
\frac{B\left(s\right)^{\frac{m}{k}}}{h\left(s\right)} \,<\, \frac{1-k}{x_3} \quad \text{for all} \quad s \in [\lambda, s_\gamma].
\end{equation}
Next, assume, by way of contradiction,
that the trajectory $T$ intersects the curve $x \,=\, \left(1-k\right)h\left(s\right)$ for some $s \in [\lambda, s_\gamma]$.
Using claim \eqref{eq:lemma21_claim} we then obtain $x_3 B\left(s\right)^{\frac{m}{k}} \,<\, x$ for the point of intersection.
But from \eqref{e7a} in Lemma~\ref{le:lemma20} it follows that $x \,<\, x_3 B\left(s\right)^{\frac{m}{k}}$ as long as $T$ stays in the region defined by
$x \,<\, \left(1-k\right)h\left(s\right)$. Using continuity this leads to a contradiction.
Hence, we conclude that the trajectory $T$ intersects $x\,=\,\left(1-k\right)\, h\left(s\right)$ next time, after $P_3$, for $s\,>\,s_\gamma$
and $T$ is inside the region determined by $x\,<\,\left(1-k\right)\, h\left(s\right)$ before it intersects $s\,=\,s_\gamma$.

To finish the proof of Lemma~\ref{le:lemma21} it remains to prove that claim \eqref{eq:lemma21_claim} holds true.
To do so we observe that differentiating $G\left(s\right)\,=\,\frac{F\left(s\right)^\frac{m}{k}}{h\left(s\right)}$ with respect to $s$,
where $F\left(s\right)$ is given by \eqref{eq:F_def},
gives
$$
G'\left(s\right)\,=\,G_*(s) \, G^*\left(s\right),
$$
where
$$
G_*(s) \, = \, \frac{ m F\left(s\right)^\frac{m}{k}}{k\left(1-s\right)^2 s \left(s+a\right)^2}
$$
and
$$
G^*\left(s\right)\,=\,2 \frac{k}{m} s^2 +\left(a \frac{k}{m} - \frac{k}{m} + 1\right)s-\lambda.
$$
As $G_*(s) > 0$ the function $G(s)$ may have max/min only if $G^*(s) = 0$. 
%
%
As $G^*(\lambda) < 0 < G^*(1)$ and $G^*(s)$ is convex
we conclude that $G\left(s\right)$ has a minimum between $\lambda$ and 1 and no other extremum.
Thus,
the maximal value of $G$ in $\lbrack \lambda ,s_\gamma\rbrack$ is either $G\left(\lambda\right)$ or $G\left(s_\gamma\right)$.
Claim \eqref{eq:lemma21_claim} now follows since
$
\frac{G\left(s\right)}{F\left(\lambda\right)^{\frac{m}{k}}} \,=\, \frac{B\left(s\right)^\frac{m}{k}}{h\left(s\right)}
$
%
and the assumptions in the lemma equals
$$
\frac{B\left(\lambda\right)^\frac{m}{k}}{h\left(\lambda\right)} \,<\, \frac{1-k}{x_3} \quad \textrm{and} \quad \frac{B\left(s_\gamma\right)^\frac{m}{k}}{h\left(s_\gamma\right)} \,<\, \frac{1-k}{x_3}.
$$
The proof of Lemma~\ref{le:lemma21} is complete. $\hfill\Box$\\

\noindent
{\it Proof of Lemma \ref{le:lemma19}}.
With notation $y = \frac{\tilde x_1}{h(\lambda)}$ it follows that
\begin{align*}
e^y\frac{\partial \bar \eta}{\partial \lambda} &= \frac{\partial}{\partial \lambda} \left( e^{\frac{\lambda}{s_\gamma}} \frac{s_\gamma + a}{1-s_\gamma} \frac{1}{a+\lambda} \right)^{\frac{m}{k}} z_i \tilde x_1\\
&+ \left( e^{\frac{\lambda}{s_\gamma}} \frac{s_\gamma + a}{1-s_\gamma} \frac{1}{a+\lambda} \right)^{\frac{m}{k}} \left( \frac{\partial z_i}{\partial \lambda} \tilde x_1
+  z_i \frac{\partial \tilde x_1}{\partial \lambda}
-z_i \tilde x_1 \frac{dy}{d\lambda}\right),
\end{align*}
in which the first term yields
\begin{align*}
&\frac{\partial}{\partial \lambda} \left( e^{\frac{\lambda}{s_\gamma}} \frac{s_\gamma + a}{1-s_\gamma} \frac{1}{a+\lambda} \right)^{\frac{m}{k}} z_i \tilde x_1\\
&= \frac{m}{k}  \left( e^{\frac{\lambda}{s_\gamma}} \frac{s_\gamma + a}{1-s_\gamma} \frac{1}{a+\lambda} \right)^{\frac{m}{k} - 1} \left( \frac{1}{s_\gamma} e^{\frac{\lambda}{s_\gamma}} \frac{s_\gamma + a}{1-s_\gamma} \frac{1}{a+\lambda}   -   e^{\frac{\lambda}{s_\gamma}} \frac{s_\gamma + a}{1-s_\gamma} \frac{1}{(a+\lambda)^2} \right) z_i \tilde x_1\\
&=\frac{m}{k}  \left( e^{\frac{\lambda}{s_\gamma}} \frac{s_\gamma + a}{1-s_\gamma} \frac{1}{a+\lambda} \right)^{\frac{m}{k} - 1}
e^{\frac{\lambda}{s_\gamma}} \frac{s_\gamma + a}{1-s_\gamma} \frac{1}{a+\lambda} \left( \frac{1}{s_\gamma}    -    \frac{1}{a+\lambda} \right) z_i \tilde x_1\\
&= \left( e^{\frac{\lambda}{s_\gamma}} \frac{s_\gamma + a}{1-s_\gamma} \frac{1}{a+\lambda} \right)^{\frac{m}{k}}
 \frac{m}{k}\left( \frac{1}{s_\gamma}    -    \frac{1}{a+\lambda} \right) z_i \tilde x_1.
\end{align*}
Therefore,
\begin{align*}
e^{\frac{\tilde x_1}{h(\lambda)}}\frac{\partial \bar \eta}{\partial \lambda} &=  \left( e^{\frac{\lambda}{s_\gamma}} \frac{s_\gamma + a}{1-s_\gamma} \frac{1}{a+\lambda} \right)^{\frac{m}{k}} \times \\
&\times\left( \frac{\partial z_i}{\partial \lambda} \tilde x_1
+  z_i \frac{\partial \tilde x_1}{\partial \lambda}
-z_i \tilde x_1 \frac{dy}{d\lambda}
+  \frac{m}{k}\left( \frac{1}{s_\gamma}    -    \frac{1}{a+ \lambda} \right) z_i \tilde x_1
\right),
\end{align*}
and so it suffices to ensure that the following expression is nonnegative:
\begin{align}\label{eq:deriv-simple}
\frac{\partial z_i}{\partial \lambda} \tilde x_1
+  z_i \frac{\partial \tilde x_1}{\partial \lambda}
+z_i \tilde x_1 \left( -\frac{dy}{d\lambda}
+  \frac{m}{k} \frac{1}{s_\gamma}    -     \frac{m}{k} \frac{1}{a+\lambda} \right).
\end{align}
With $\tilde x_1$ as in \eqref{eq:etasimp1} and \eqref{eq:x1simp1} we have, where $s_0$ is either $0.475$ or $0.8$, that
\begin{align}\label{eq:deriv:1}
  \frac{\partial \tilde x_1}{\partial \lambda} = m \ln\frac{\lambda}{s_0} 
\end{align}
Further, $Y = ye^{-y}$ and $y = \tilde x_1 / h(\lambda)$ so that
\begin{align}\label{eq:deriv:y}
\frac{\partial y}{\partial \lambda} = - \frac{\tilde x_1}{h(\lambda)^2} + \frac{\partial \tilde x_1}{\partial \lambda} \frac{1}{h(\lambda)}
= - \frac{\tilde x_1}{h(\lambda)^2} +  \frac{m}{h(\lambda)} \ln\frac{\lambda}{s_0}
\end{align}
and so, for $i = 0,2$, we have
\begin{align}\label{eq:deriv:2}
\frac{\partial z_i}{\partial \lambda} &=  \frac{d z_i}{d Y}\frac{d Y}{d y} \frac{\partial y}{\partial \lambda}
= (1-y)e^{-y}\frac{d z_i}{d Y}  \frac{\partial y}{\partial \lambda}\notag \\
&= (y-1)e^{-y}\frac{d z_i}{d Y}  \left(- \frac{m}{h(\lambda)}  \ln\frac{\lambda}{s_0}  + \frac{\tilde x_1}{h(\lambda)^2}\right).
\end{align}
Further, as $z_0 = (1-(e-1) Y)^{-1}$ we have
\begin{align*}
\frac{d z_0}{d Y} = \frac {e-1}{(1-(e-1) Y)^2} > 0,
\end{align*}
and by construction (Recall the proof of Lemma \ref{le:julelemma})
$$
\frac{f_2(z)}{z} = Y \quad \text{in which} \quad f_2(z) = \frac{z-1}{c_2 z + d_2}.
$$
Differentiation of $ c_2 z^2 Y -(1- d_2) z  + 1 = 0$ yields
$$
\frac{d z_2}{d Y} = -\frac{ c2 z^2 + d_2 z}{2 c_2 Y z - (1 - d_2 Y)} > 0
$$
where the the last step follows since $c_2 = \frac{1}{e}, d_2 = e-2$, $Y \in (0,\frac{1}{e})$ and $z \in (1,e)$.
Thus
\begin{align}\label{eq:deriv:z0z2>0}
\frac{d z_i}{d Y} > 0 \quad \text{for} \quad i = 1,2.
\end{align}
Using the derivatives in \eqref{eq:deriv:1}, \eqref{eq:deriv:y} and \eqref{eq:deriv:2} we realize that \eqref{eq:deriv-simple} is equivalent with
\begin{align*}
(y-1)e^{-y}\frac{d z_i}{d Y} \left(- \frac{m}{h(\lambda)}  \ln\frac{\lambda}{s_0}  + \frac{\tilde x_1}{h(\lambda)^2}\right) \tilde x_1
+  z_i m \ln \frac{\lambda}{s_0} \\
+z_i \tilde x_1 \left( \frac{\tilde x_1}{h(\lambda)^2} - \frac{m}{h(\lambda)} \ln \frac{\lambda}{s_0}
+  \frac{m}{k} \frac{1}{s_\gamma}    -     \frac{m}{k} \frac{1}{a+\lambda} \right).
\end{align*}
By recalling \eqref{eq:deriv:z0z2>0} and excluding the obviously positive terms we conclude that it suffices to ensure that
\begin{align*}
  \ln \frac{\lambda}{s_0}
  - \frac{\tilde x_1}{h(\lambda)} \ln \frac{\lambda}{s_0}
    -     \frac{1}{k} \frac{\tilde x_1}{h(\lambda)}  \geq 0,
\end{align*}
which follows by
\begin{align*}
\kappa\, \tilde x_1 > a + \lambda > h(\lambda) \quad \text{and} \quad  (1-\kappa) k \ln  \frac{s_0}{\lambda} > 1
\end{align*}
for any $\kappa \in (0,1)$.

\bigskip

We now turn to the monotonicity in $a$.
It follows that
\begin{align*}
e^{\frac{\tilde x_1}{h(\lambda)}}\frac{\partial \bar \eta}{\partial a} &= \frac{\partial}{\partial a} \left( e^{\frac{\lambda}{s_\gamma}} \frac{s_\gamma + a}{1-s_\gamma} \frac{1}{a+\lambda} \right)^{\frac{m}{k}} z_i \tilde x_1\\
&+ \left( e^{\frac{\lambda}{s_\gamma}} \frac{s_\gamma + a}{1-s_\gamma} \frac{1}{a+\lambda} \right)^{\frac{m}{k}} \left( \frac{\partial z_i}{\partial a} \tilde x_1
+  z_i \frac{\partial \tilde x_1}{\partial a}
- z_i \tilde x_1  \frac{\partial y}{\partial a}\right),
\end{align*}
in which the first term yields
\begin{align*}
&\frac{\partial}{\partial a} \left( e^{\frac{\lambda}{s_\gamma}} \frac{s_\gamma + a}{1-s_\gamma} \frac{1}{a+\lambda} \right)^{\frac{m}{k}} z_i \tilde x_1\\
&= \frac{m}{k}  \left( e^{\frac{\lambda}{s_\gamma}} \frac{s_\gamma + a}{1-s_\gamma} \frac{1}{a+\lambda} \right)^{\frac{m}{k} - 1}
\left( e^{\frac{\lambda}{s_\gamma}} \frac{1}{1-s_\gamma} \frac{1}{a+\lambda}\left(1   -   \frac{s_\gamma + a}{a+\lambda}\right) \right) z_i \tilde x_1\\
&= \left( e^{\frac{\lambda}{s_\gamma}} \frac{s_\gamma + a}{1-s_\gamma} \frac{1}{a+\lambda} \right)^{\frac{m}{k}}
 \frac{m}{k}\left( \frac{1}{s_\gamma+a}    -    \frac{1}{a+\lambda} \right) z_i \tilde x_1.
\end{align*}
Therefore,
\begin{align*}
e^{\frac{\tilde x_1}{h(\lambda)}}\frac{\partial \bar \eta}{\partial a} &=  \left( e^{\frac{\lambda}{s_\gamma}} \frac{s_\gamma + a}{1-s_\gamma} \frac{1}{a+\lambda} \right)^{\frac{m}{k}} \times \\
&\times\left( \frac{\partial z_i}{\partial a} \tilde x_1
+  z_i \frac{\partial \tilde x_1}{\partial a}
- z_i \tilde x_1 \frac{\partial y}{\partial a}
+  \frac{m}{k}\left( \frac{1}{s_\gamma + a}    -    \frac{1}{a+\lambda} \right) z_i \tilde x_1
\right),
\end{align*}
and so it suffices to ensure that the following expression is nonnegative:
\begin{align}\label{eq:deriv-simple:a}
\frac{\partial z_i}{\partial a} \tilde x_1
+  z_i \frac{\partial \tilde x_1}{\partial a}
+z_i \tilde x_1 \left( - \frac{\partial y}{\partial a}
+  \frac{m}{k} \frac{1}{s_\gamma + a}    -     \frac{m}{k} \frac{1}{a+\lambda} \right).
\end{align}
Recall $y = \tilde x_1 / h(\lambda)$ and so 
\begin{align}\label{eq:derivya}
\frac{\partial y}{\partial a} = - \frac{\tilde x_1}{h(\lambda)^2} + \frac{\partial \tilde x_1}{\partial a} \frac{1}{h(\lambda)}
%
%
\end{align}
and since $Y = ye^{-y}$ we also conclude,
for $i = 0,2$, that
\begin{align}\label{eq:deriv:2a}
\frac{\partial z_i}{\partial a} &=  \frac{d z_i}{d Y}\frac{d Y}{d y} \frac{\partial y}{\partial a}
= (1-y)e^{-y}\frac{d z_i}{d Y}  \frac{\partial y}{\partial a}
%
\end{align}
Inserting \eqref{eq:derivya} and \eqref{eq:deriv:2a} into \eqref{eq:deriv-simple:a} gives
\begin{align}\label{eq:deriv-simple:aa}
&\tilde x_1 (y-1)e^{-y}\frac{d z_i}{d Y}  \left(- \frac{\partial \tilde x_1}{\partial a} \frac{1}{h(\lambda)}  + \frac{\tilde x_1}{h(\lambda)^2}\right)
+  z_i \frac{\partial \tilde x_1}{\partial a}\notag\\
& + z_i \tilde x_1 \left(- \frac{\partial \tilde x_1}{\partial a} \frac{1}{h(\lambda)}  + \frac{\tilde x_1}{h(\lambda)^2}
+  \frac{m}{k} \frac{1}{s_\gamma + a}    -     \frac{m}{k} \frac{1}{a+\lambda} \right).
\end{align}

We begin with the case $m < 0.3$. Then from \eqref{eq:x1simp1} we have $\frac{\partial \tilde x_1}{\partial a} = 0$ and hence \eqref{eq:deriv-simple:aa} simplifies to
\begin{align*}
\tilde x_1 (y-1)e^{-y}\frac{d z_i}{d Y} \frac{\tilde x_1}{h(\lambda)^2}
+  z_i \frac{\partial \tilde x_1}{\partial a}
 z_i \tilde x_1 \left( \frac{\tilde x_1}{h(\lambda)^2}
+  \frac{m}{k} \frac{1}{s_\gamma + a}    -     \frac{m}{k} \frac{1}{a+\lambda} \right).
\end{align*}
which is nonnegative (recall $y > 1$ and \eqref{eq:deriv:z0z2>0}) if
\begin{align*}
\frac{m}{k}  \leq \frac{\tilde x_1}{a + \lambda}. 
\end{align*}
%
Plugging in $k \geq 1/3$ and $m \leq 0.3$ we see that it suffices to have $\tilde x_1 \leq a+\lambda$ which holds by assumption.

When $m \geq 0.3$ it follows from \eqref{eq:x1simp1} that $\frac{\partial \tilde x_1}{\partial a} = 0.2$ and hence \eqref{eq:deriv-simple:aa} simplifies to
\begin{align*}
&\tilde x_1 (y-1)e^{-y}\frac{d z_i}{d Y}  \left(- \frac{1}{5 h(\lambda)}  + \frac{\tilde x_1}{h(\lambda)^2}\right)
+  z_i \frac{\partial \tilde x_1}{\partial a}\notag\\
& + z_i \tilde x_1 \left(- \frac{1}{5 h(\lambda)}  + \frac{\tilde x_1}{h(\lambda)^2}
+  \frac{m}{k} \frac{1}{s_\gamma + a}    -     \frac{m}{k} \frac{1}{a+\lambda} \right)
\end{align*}
which is nonnegative if
$$
\frac{1}{5} + \frac{m}{k} \frac{h(\lambda)}{a+\lambda} \leq \frac{\tilde x_1}{h(\lambda)} = \frac{h(0.8) + m\, (0.8 -\lambda\, (1-\ln \lambda +\ln 0.8 ))}{h(\lambda)}
$$
which is implied by
$$
\frac{1}{5} + \frac{m}{k}  \leq  \frac{\frac{1}{5}(a + \lambda) + 0.6 m}{a + \lambda}
$$
and, thanks to assumptions $k\geq 1/3$ and $a + \lambda \leq 0.2$, we conclude the proof. $\hfill \Box$


\subsection{Numerical results in Region 4}
\label{sec:numerics_reg4}

\begin{figure}[h]
\begin{center}
\includegraphics[scale = 0.5]{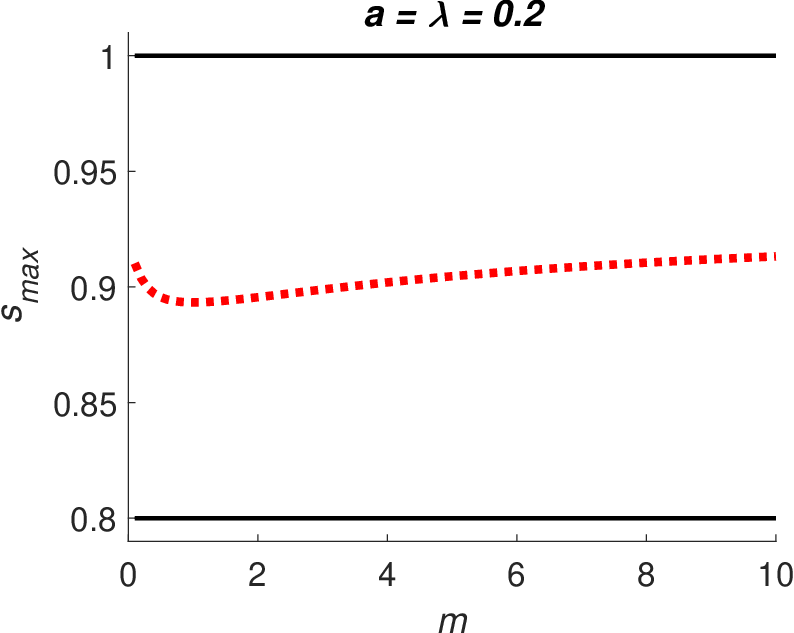}\hspace{0.5cm}
\includegraphics[scale = 0.5]{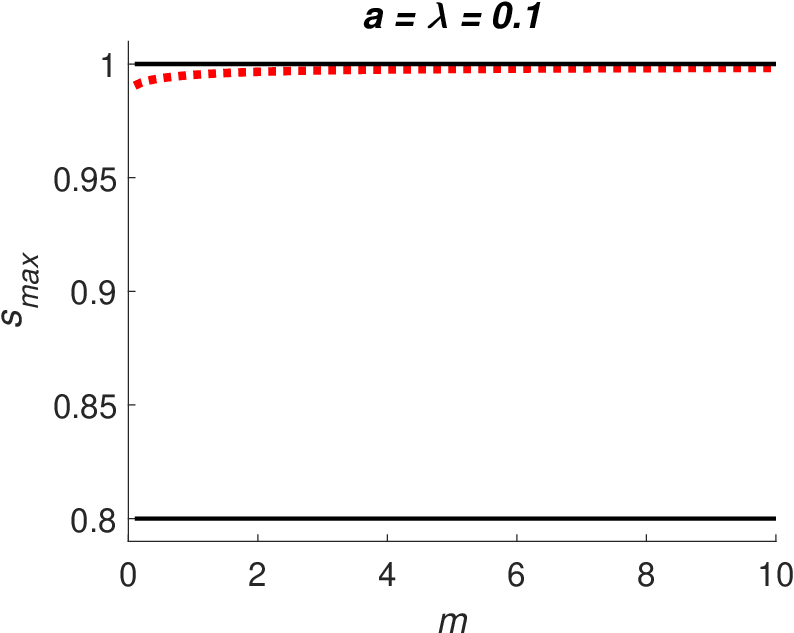}\vspace{1cm}
\includegraphics[scale = 0.5]{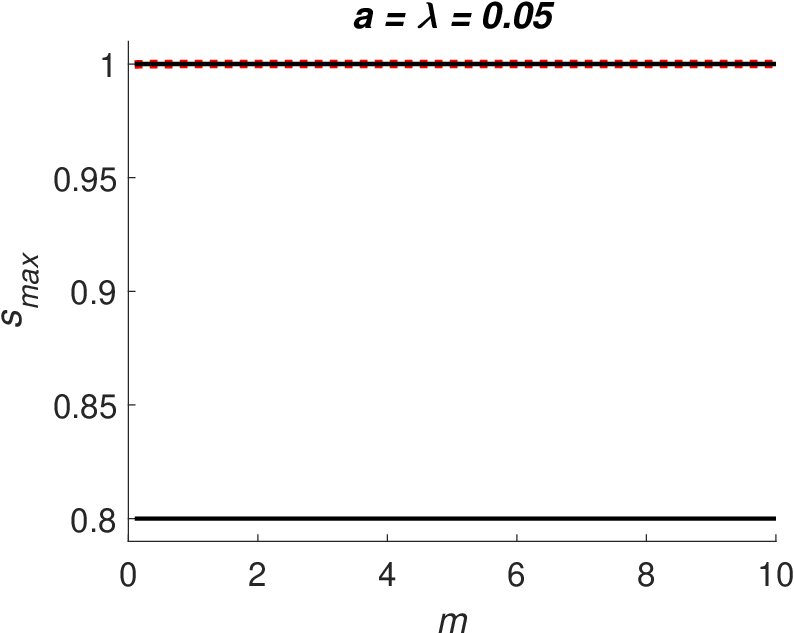}\hspace{0.5cm}
\includegraphics[scale = 0.5]{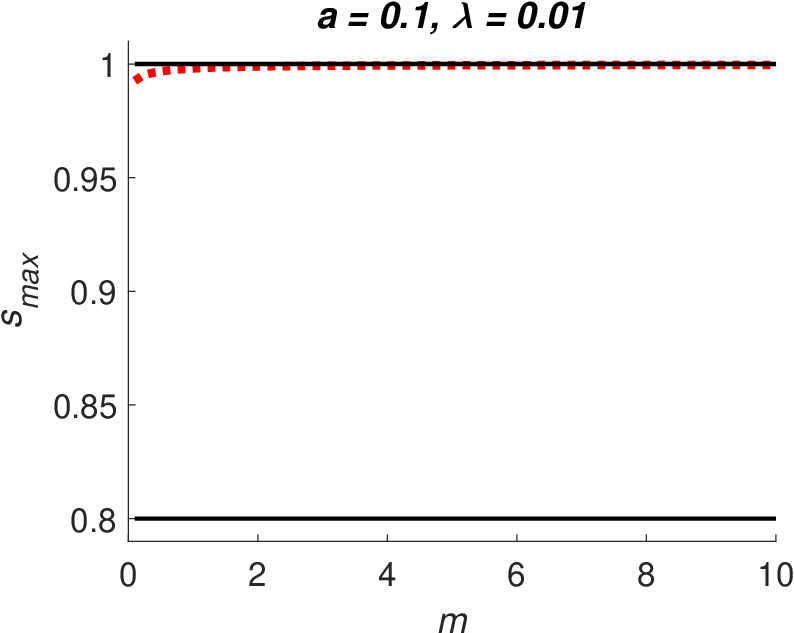}
\caption{
The estimates 0.8 and 1 (black, solid) together with $s_{max}$ on (simulations of) the stable limit cycle.}
\label{fig:smax}
\end{center}
\end{figure}

As in subsections \ref{sec:numerics_reg1} and \ref{sec:numerics_reg23} we end the section with numerical simulations.
Figure~\ref{fig:smax} shows $s_{max}$ on simulations of the stable limit cycle together with the estimates 0.8 and 1 as
functions of $m$ for the same parameter values as in Fig. \ref{fig:xmax}.


\begin{thebibliography}{99}

\bibitem{BHW83}
Butler, G. J.,  Hsu, S. B.,  Waltman, P.:
{Coexistence of competing predators in a chemostat}.
Journal of Mathematical Biology, 17.2, 1983, 133-151.

\bibitem{C81}
Cheng, K.-S.:
{Uniqueness of a limit cycle for a predator-prey system}.
SIAM Journal on Mathematical Analysis, 12.4, 1981, 541-548.

\bibitem{Angelis75}
De Angelis, D. L.
{Estimates of predator-prey limit cycles.}
Bulletin of Mathematical Biology, 37, 1975, 291-299.

%
\bibitem{EOS96}
Eirola, T., Osipov, A.V., S\"oderbacka, G.:
{Chaotic regimes in a dynamical system of the type many predators one prey}.
Research reports A, 386, Helsinki University of Technology, 1996.
%
%
%
%
%

\bibitem{H00}
Hasik K.
{Uniqueness of limit cycle in the predator–prey system with symmetric prey isocline}.
Mathematical biosciences, 164.2, 2000, 203-215.


\bibitem{HS09}
Hsu, S.-B., Shi, J.:
{Relaxation oscillator profile of limit cycle in predator-prey system}.
Discrete and continuous dynamical systems B, 11, no. 4, 2009, 893-911.

\bibitem{HW20}
Hsu, T.-H., Wolkowicz, G.S.K.
{A criterion for the existence of relaxation oscillations with applications to predator-prey systems and an epidemic model}.
Discrete and Continuous Dynamical Systems B, 25.4, 2020, 1257-1277.

\bibitem{H88}
Huang, X.-C.:
{Uniqueness of limit cycles of generalised Li\'enard systems and predator-prey systems}.
Journal of Physics A: Mathematical and General, 21.13 (1988): L685.

\bibitem{HM89}
Huang, X.-C., Merrill, S. J.:
{Conditions for uniqueness of limit cycles in general predator-prey systems}.
Mathematical biosciences, 96.1 (1989): 47-60.


\bibitem{KF88}
Kuang, Y., Freedman, H. I.:
{Uniqueness of limit cycles in Gause-type models of predator-prey systems}.
Mathematical Biosciences, 88.1 (1988): 67-84.

\bibitem{ns}
Lundstr\"om N. L. P., S\"oderbacka G., Estimates of size of cycle in a predator-prey system,
Differential Equations and Dynamical Systems, 30, 2022, 131-159.

\bibitem{LS23}
Lundstr\"om N. L. P., S\"oderbacka G. J., Analytical approximations of Lotka-Volterra integrals,
arXiv preprint arXiv:2311.03069, 2023.





 \bibitem{osipovAlushta}
Osipov A. V.,  S{\"o}derbacka G.,
      Extinction and coexistence of predators,
     { Dinamicheskie Sistemy},
      {6}, 1, 2016, 55-64.

 \bibitem{osipovIJBC}
Osipov A. V.,   S{\"o}derbacka G.,
    Poincar\'{e}  map construction for some classic two predators - one prey systems,
      { Internat. J. Bifur. Chaos Appl. Sci. Engrg},
       {27}, 2017, no 8, 1750116, 9.

 \bibitem{osipoveuler}
Osipov A. V.,   S{\"o}derbacka G.,
Review of results on a system of type many predators-one prey,
{Nonlinear Systems and Their Remarkable Mathematical Structures}, 2018, 520-540.
CRC Press.

\bibitem{RMK93}
Rinaldi, S., Muratori, S., Kuznetsov, Y.:
{Multiple attractors, catastrophes and chaos in seasonally perturbed predator-prey communities}.
Bulletin of mathematical Biology, 55.1, 1993, 15-35.

\bibitem{R71}
Rosenzweig, M. L.:
Paradox of enrichment: destabilization of exploitation ecosystems in ecological time.
Science, 171.3969, 1971, 385-387.

\bibitem{S23}
S\"oderbacka G. J.,
{Model map and multistability for a two predator-one prey system}.
Differential Equations and Control Processes, 1, 2023, 12-23

\bibitem{RM63}
Rosenzweig, M. L., MacArthur R.:
{Graphical representation and stability conditions of predator-prey interaction}.
American Naturalist, 97, 1963, 209-223.


%
%

%
%
%
\end{thebibliography}
\end{document}